\documentclass[leqno]{amsart}
\usepackage{egastyle}
\usepackage{mymacros}
\usepackage[all,tips]{xy}
\usepackage{url}
\usepackage{bm}
\usepackage{verbatim}

\CompileMatrices

\def\cA{{\mathcal A}}

\def\cX{{\mathcal X}}

\newcommand{\PP}{\mathbf{P}}

\newcommand{\ZZ}{\mathbf{Z}}
\newcommand{\HH}{\mathbf{H}}

\renewcommand{\setminus}{\smallsetminus}

\newcommand{\an}{{\rm an}}

\newcommand{\berk}{{\rm an}}

\renewcommand{\Div}{{\rm div}}

\newcommand\pu[1]{\{\kern -2pt\{#1\}\kern -2pt\}}

\DeclareMathOperator{\DV}{DV}
\DeclareMathOperator{\val}{val}

\DeclareMathOperator{\trop}{trop}

\DeclareMathOperator{\diam}{diam}

\newcommand\A{\bA}
\newcommand\B{\bB}
\newcommand\G{\bG}

\title{On the structure of nonarchimedean analytic curves}

\author{Matthew Baker} 
\email{mbaker@math.gatech.edu}
\address{School of Mathematics, Georgia Institute of Technology, Atlanta GA 30332-0160, USA}
\author{Sam Payne}
\email{sam.payne@yale.edu}
\address{Mathematics Department, Yale University, New Haven, CT 06511}
\author{Joseph Rabinoff}
\email{rabinoff@math.harvard.edu}
\address{Department of Mathematics, Harvard University, Cambridge, MA 02138}

\begin{document}

\begin{abstract}
  Let $K$ be an algebraically closed, complete nonarchimedean field and
  let $X$ be a smooth $K$-curve.  In this paper we elaborate on several
  aspects of the structure of the Berkovich analytic space
  $X^\an$.  We define semistable vertex sets of $X^\an$ and their
  associated skeleta, which are essentially finite metric graphs embedded
  in $X^\an$.  We prove a folklore theorem which states that semistable
  vertex sets of $X$ are in natural bijective correspondence with
  semistable models of $X$, thus showing that our notion of skeleton
  coincides with the standard definition of
  Berkovich~\cite{berkovich:analytic_geometry}.  We use the skeletal
  theory to define a 
  canonical metric on $\HH(X^\an) \coloneq X^\an\setminus X(K)$, and we
  give a proof of Thuillier's nonarchimedean Poincar\'e-Lelong formula in
  this language using results of Bosch and L\"utkebohmert.
\end{abstract}

\maketitle

\vspace{-20 pt}

\section{Introduction}

Throughout this paper we let $K$ denote an algebraically closed field
which is complete with respect to a nontrivial, nonarchimedean valuation
$\val:K\to\R\cup\{\infty\}$.  Let $R$ be the valuation
ring of $K$, $\fm_R$ its maximal ideal, 
and $k$ its residue field.  In this situation, $R$ is not Noetherian and $k$
is algebraically closed.  We let $|\cdot| = \exp(-\val(\cdot))$ be an
associated absolute value and let $G = \val(K^\times)\subset\R$ be the
value group.

\paragraph
Let $X$ be a smooth, proper, connected algebraic $K$-curve and let $X^\an$
be its analytification in the sense of
Berkovich~\cite{berkovich:analytic_geometry}.  
The purpose of this note is to elaborate on the following aspects of the
structure of 
$X^\an$:
\begin{enumerate}
\item We define semistable vertex sets of $X$ and their associated
  skeleta, which are finite metric graphs contained in $X^\an$.
\item We make explicit the bijective correspondence between the semistable
  vertex sets of $X$ and the semistable models of $X$.
\item We show that `most' of $X$ can be covered by  skeleta as
  above, and  use this fact to define a canonical metric on 
  $\HH(X^\an) \coloneq X^\an\setminus X(K)$ such that the resulting metric
  space is locally 
  modeled on an $\R$-tree. 
\item We use results of Bosch and L\"utkebohmert to prove
  Thuillier's nonarchimedean analogue of the Poincar\'e-Lelong formula.
  This includes the fact that the valuation of a nonzero rational function on $X$
  is a piecewise affine function on $\HH(X^\an)$, which is moreover
  \emph{harmonic} in the sense that the sum of the outgoing slopes at any
  point is zero.  
\end{enumerate}
The above results, while very useful and in large part
well-known to experts, are often difficult or impossible to extract from
the literature, although Ducros~\cite{ducros:courbes_analytiques} is
preparing a book on the subject.
In~\cite{bpr:trop_curves} we apply these ideas to study the
relationship between the analytification of a curve $X$ and its
tropicalization with respect to a rational map to a torus.  In particular, we
study the metric aspects of the tropicalization map and the
relationship between skeleta and tropicalizations of $X$.

\paragraph[Skeleta and semistable models]
The theory of skeleta goes back to
Berkovich~\cite{berkovich:analytic_geometry,berkovich:locallycontractible1},
and is elaborated 
somewhat in the case of curves in Thuillier's (unpublished)
thesis~\cite{thuillier:thesis}. 
In Berkovich's approach, a skeleton is a subset of $X^\an$ which is
associated to a semistable model of $X$.  In contrast, we define the
skeleton in terms of 
a \emph{semistable vertex set} of $X$
(called a \emph{triangulation} in~\cite{ducros:triangulations}; see also
Temkin's lecture notes~\cite{temkin:berkovich_spaces}), which is
a finite set $V$ of type-$2$ 
points of $X^\an$ such that $X^\an\setminus V$ is a disjoint union of
open balls and finitely many open annuli.  This approach has the
advantage that it only makes reference to the analytic space $X^\an$ and
is not conceptually tied to the semistable reduction theory of $X$, thus
making certain constructions more natural.

If $V\subset X^\an$ is a semistable vertex set, the connected component
decomposition 
\[ X^\an \setminus V = \Djunion \B(1)_+ \djunion 
\Djunion_{i=1}^r \bS(a_i)_+ \]
is called a \emph{semistable decomposition} of $X$; here 
$\B(1)_+$ is the open unit ball and for $a_i\in\fm_R\setminus\{0\}$,
we define $\bS(a_i)_+$ to be the open annulus of inner radius $|a_i|$ and
outer radius $1$.  A semistable decomposition of a nonarchimedean curve is
somewhat analogous to a pair-of-pants decomposition of a Riemann surface.
The annulus $\bS(a_i)_+$ has a canonical closed subset
$\Sigma(\bS(a_i)_+)$, called its \emph{skeleton}, which is 
identified with the open interval $(0,\val(a))$.  The
\emph{skeleton} of $X$ associated to $V$ is the set
\[ \Sigma(X,V) \coloneq V \cup \bigcup_{i=1}^r \Sigma(\bS(a_i)_+). \]
We show that $\Sigma(X,V)$ is naturally homeomorphic to a finite graph,
with vertices 
$V$ and open edges $\{\Sigma(\bS(a_i)_+)\}_{i=1}^r$.  Declaring 
the length of $\Sigma(\bS(a_i)_+)$ to be $\val(a_i)$ (the logarithmic
modulus of $\bS(a_i)_+$) makes $\Sigma(X,V)$ into a metric graph.
There is a deformation retraction $X^\an\to\Sigma(X,V)$, so
$\Sigma(X,V)$ is connected.
We give a relatively compelete account of these skeleta in
\S\ref{sec:skeleta.curves}. 

\paragraph
Now let $\cX$ be a semistable $R$-model of $X$.  There is a 
\emph{reduction map} $\red: X^\an\to\cX_k$, defined as follows.  Points of
$X^\an$ correspond in a natural way to equivalence classes of maps 
$\Spec(L)\to X$, where $L/K$ is a \emph{valued} field extension.  By the
valuative criterion of properness, a point $x:\Spec(L)\to X$ extends in a
unique way to a map $\Spec(\sO_L)\to\cX$ from the valuation ring of $L$; the
reduction $\red(x)$ is defined to be the image of the closed point.
The reduction map is anti-continuous, in that the inverse image of a
Zariski-open set in $\cX_k$ is a closed subset of $X^\an$ (and vice
versa).  The fibers of $\red$ are called \emph{formal fibers}.

It follows from a theorem of Berkovich that for a generic point
$\bar x\in\cX_k$, the formal fiber $\red\inv(\bar x)$ consists of a single
type-$2$ point.  If $\bar x\in\cX_k$ is a closed point then by a
theorem of Bosch and L\"utkebohmert, the formal fiber $\red\inv(\bar x)$
is isomorphic to an open ball or an open annulus if $\bar x$ is a smooth
point or a node, respectively.  It follows from this and the
anti-continuity of $\red$ that the set $V_\cX$ of points
of $X^\an$ reducing to generic points of $\cX_k$ is a semistable vertex
set, and that the decomposition
\[ X^\an\setminus V = \Djunion_{\bar x\in\cX(k)} \red\inv(\bar x) \]
of $X^\an$ into its formal fibers is a semistable decomposition.  The
associated skeleton $\Sigma(X,V_\cX)$ is the incidence graph of 
the irreducible components of $\cX_k$.

We will give a proof of the following folklore theorem which says that the
associaton $\cX\mapsto V_\cX$ is bijective, thus unifying the two
notions of skeleta.
See Theorem~\ref{thm:construct.models} and Proposition~\ref{prop:skeleta.agree}.

\begin{thm*}
  The association $\cX\mapsto V_\cX$ is a bijective correspondence from
  the set of semistable models of $X$ to the set of semistable vertex 
  sets of $X$.  Moreover, there exists a morphism  of semistable models
  $\cX\to\cX'$ if and only if $V_{\cX'}\subset V_{\cX}$.
\end{thm*}

\paragraph
In \S\ref{sec:metric.structure} we turn our attention to the metric
nature of the analytic curve $X^\an$.  This is worked out `by hand' in
the case $X = \PP^1$ by Baker and Rumely~\cite{baker_rumely:book}, but does
not otherwise explicitly appear in the literature.  We will
prove that the metric graph structures on each skeleton of $X^\an$
are compatible, and that the resulting metric on their union extends by continuity
to a unique metric on $\HH(X^\an) = X^\an\setminus X(K)$.  The resulting metric
space is locally modeled on an $\R$-tree.  At this point it is
straightforward to give a proof of Thuillier's 
Poincar\'e-Lelong formula~\cite[Proposition~3.3.15]{thuillier:thesis}
using classical machinery of Bosch and L\"utkebohmert.  We will do so
without developing the harmonic analysis necessary to give the statement
of Thuillier's theorem in terms of the nonarchimedean $dd^c$ operator; for
this reason we call the theorem the Slope Formula.

In the statement of the Slope Formula, we say that a function
$F:\HH(X^\an)\to\R$ is \emph{piecewise affine with integer slopes}
provided that, for every isometric embedding
$\alpha:[a,b]\inject\HH(X^\an)$, the composition $F\circ\alpha$ is
a piecewise affine function with integer slopes on the interval $[a,b]$.
We define the \emph{set of tangent directions $T_x$} at a point $x\in\HH(X^\an)$ to be the
set of germs of isometric embeddings $\alpha:[a,b]\inject X^\an$ such that
$\alpha(a) = x$, and we define the \emph{outgoing slope of $F$} in the
tangent direction $v$ represented by $\alpha$ to be the right-hand
derivative $d_vF(x)$ of $F\circ\alpha$ at $a$.  One can extend these
definitions to closed points $x\in X(K)$, although in our later
formulation of the Slope Formula we avoid this issue.

\begin{thm*}[Slope Formula]
  Let $f\in K(X)^\times$ and let 
  $F = -\log|f|:X^\an\to\R\cup\{\pm\infty\}$.
  \begin{enumerate}
  \item $F$ is piecewise affine with integer slopes on $\HH(X^\an)$.
  \item For $x\in\HH(X^\an)$ we have $d_v F(x) = 0$ for almost all $v\in
    T_x$, and 
    \[ \sum_{v\in T_x} d_v F(x) = 0. \]
    In other words, $F$ is \emph{harmonic}.
  \item For $x\in X(K)$ there is a unique tangent direction $v\in T_x$ and
    $d_v F(x) = \ord_x(f)$.
  \end{enumerate}
\end{thm*}

\smallskip
See Theorem~\ref{thm:PL} for a more precise statement and proof.

\paragraph[Acknowledgements]
\thanks{The authors would like to thank Vladimir Berkovich, Antoine
  Ducros, Walter Gubler, and Michael Temkin for helpful and illuminating
  discussions.}

\section{The skeleton of a generalized annulus}

In this section we prove some preliminary facts about the building blocks of
analytic curves, namely, open balls and open
annuli.  We also study punctured open balls in order to treat marked (or
punctured) curves and their skeleta in \S\ref{sec:skeleta.curves}.

\paragraph[Some analytic domains in $\A^1$]
Define the extended tropicalization map, or valuation map,
\[ \trop:\sM(K[T]) = \A^{1,\an}\to\R\cup\{\infty\} \sptxt{by}
\trop(\|\cdot\|) = -\log(\|T\|); \]
here $\sM(\,\cdot\,)$ denotes the Berkovich spectrum.
Clearly $\G_m^\an = \trop\inv(\R)$.
We use $\trop$ to define several analytic domains in $\A^{1,\an}$:
\begin{bullets}
\item For $a\in K^\times$ the 
  \emph{standard closed ball of radius $|a|$} is
  $\B(a) = \trop\inv([\val(a),\infty])$. 
  This is a polyhedral domain 
  in the sense of~\cite{jdr:trop_ps}; more precisely, it is the affinoid
  domain with ring of analytic functions 
  \[ K\angles{a\inv t} = 
  \left\{ \sum_{n=0}^\infty a_n t^n~:~ |a_n|\cdot|a|^n\to 0 
  \text{ as } n\to\infty\right\}. \]
  The supremum norm is given by
  \[ \left|\sum_{n=0}^\infty a_n t^n\right|_{\sup}
  = \max \big\{ |a_n|\cdot|a|^n ~:~ n\geq 0 \big\} \]
  and the canonical reduction is the polynomial ring $k[\tau]$, where
  $\tau$ is the residue of $a\inv t$.

\item For $a\in K^\times$ the
  \emph{standard open ball of radius $|a|$} is
  $\B(a)_+ = \trop\inv((\val(a),\infty])$.  This is an open analytic
  domain which can be expressed as an increasing union of standard closed
  balls.

\item For $a,b\in K^\times$ with $|a| \leq |b|$ the 
  \emph{standard closed annulus of inner radius $|a|$ and outer radius $|b|$} is 
  $\bS(a,b) = \trop\inv([\val(b),\val(a)])$.  This is a polytopal
  domain in $\G_m^\an$ ~\cite{gubler:tropical,jdr:trop_ps}; it is therefore
  an affinoid space whose ring of analytic functions is 
  \[ K\angles{a t\inv, b\inv t}
  = \left\{ \sum_{n=-\infty}^\infty a_n t^n ~:~
    |a_n|\cdot |a|^n \to 0\text{ as } n\to +\infty,~
    |a_n|\cdot |b|^n \to 0\text{ as } n\to -\infty \right\}. \]
  The supremum norm is given by
  \[ \left|\sum_{n=-\infty}^\infty a_n t^n\right|_{\sup}
  = \max \big\{ |a_n|\cdot |a|^n, |a_n|\cdot |b|^n~:~n\in\Z \big\} \]
  and the canonical reduction is $k[\sigma,\tau]/(\sigma\tau - \bar{a/b})$,
  where $\sigma$ (resp.\ $\tau$) is the residue of $a t\inv$ 
  (resp.\ $b\inv t$) and $\bar{a/b}\in k$ is the residue of $a/b$.  The canonical
  reduction is an integral domain if and only if $|a|=|b|$, in which case the
  supremum norm is multiplicative.  The \emph{(logarithmic) modulus} of
  $\bS(a,b)$ is by definition $\val(a) - \val(b)$.

\item In the above situation, if $|a| \leq 1$ and $|b| = 1$ we write
  $\bS(a) \coloneq \bS(a,1) = \trop\inv([0,\val(a)])$.  In this case
  \[ K\angles{a t\inv, t} \cong K\angles{s,t}/(st - a). \]

\item For $a,b\in K^\times$ with $|a| < |b|$ the 
  \emph{standard open annulus of inner radius $|a|$ and outer radius $|b|$} is 
  $\bS(a,b)_+ = \trop\inv((\val(b),\val(a)))$.  This is an open analytic
  domain which can be expressed as an increasing union of standard closed
  annuli.  The \emph{(logarithmic) modulus} of
  $\bS(a,b)_+$ is by definition $\val(a) - \val(b)$.  As above we write 
  $\bS(a)_+ \coloneq \bS(a,1)_+ = \trop\inv((0,\val(a)))$.

\item For $a\in K^\times$ the
  \emph{standard punctured open ball of radius $|a|$} is
  $\bS(0,a)_+ = \trop\inv((\val(a),\infty))$, and the 
  \emph{standard punctured open ball of radius $|a|\inv$ around $\infty$}
  is $\bS(a,\infty)_+ = \trop\inv((-\infty,\val(a)))$.
  These are open analytic domains which can be written as an increasing
  union of standard closed annuli.  By convention we define the \emph{modulus} of 
  $\bS(0,a)_+$ and $\bS(a,\infty)_+$ to be infinity.  
  We write $\bS(0)_+ = \bS(0,1)_+$.

\end{bullets}

\smallskip
Note that if $A$ is any of the above analytic domains in $\A^{1,\an}$
then $A = \trop\inv(\trop(A))$.
By a \emph{standard generalized annulus} we will mean a standard closed
annulus, a standard open annulus, or a standard punctured open ball, and
by a \emph{standard generalized open annulus} we will mean a standard open
annulus or a standard punctured open ball.
Note that by scaling we have isomorphisms 
\[ \B(a)\cong\B(1) \quad \B(a)_+\cong\B(1)_+ \quad
\bS(a,b)\cong\bS(ab\inv) \quad 
\bS(a,b)_+\cong\bS(ab\inv)_+ \quad
\bS(0,a)_+\cong\bS(0)_+ \]
and taking $t\mapsto t\inv$ yields $\bS(1,\infty)_+\cong\bS(0,1)_+$.

Morphisms of standard closed annuli have the following structure:

\begin{prop} \label{prop:annulus.morphism}
  Let $a\in R\setminus\{0\}$.
  \begin{enumerate}
  \item The units in $K\angles{at\inv, t}$ are the functions of the form
    \begin{equation} \label{eq:annulus.units} 
      f(t) = \alpha\,t^d (1 + g(t)) 
    \end{equation}
    where $\alpha\in K^\times$, $d\in \Z$, and $|g|_{\sup} < 1$.

  \item Let $f(t)$ be a unit as in~\eqref{eq:annulus.units} with 
    $d > 0$ (resp.\ $d < 0$).  The induced morphism
    $\phi:\bS(a)\to\bG_m^\an$ factors through a finite flat morphism 
    $\bS(a)\to \bS(\alpha a^d, \alpha)$ 
    (resp.\ $\bS(a)\to \bS(\alpha, \alpha a^d)$) of degree $|d|$.

  \item Let $f(t)$ be a unit as in~\eqref{eq:annulus.units} with $d=0$.
    The induced morphism $\phi:\bS(a)\to\bG_m^\an$ factors through a
    morphism $\bS(a)\to\bS(\alpha,\alpha)$ which is not finite.
  \end{enumerate}
\end{prop}

\pf The first assertion is proved 
in~\cite[Lemme~2.2.1]{thuillier:thesis} by considering the Newton polygon of
$f(t)$.  To prove~(2) we easily reduce to the case $\alpha = 1$ and 
$d > 0$.  Since $|f|_{\sup} = 1$ and $|f\inv|_{\sup} = |a|^{-d}$ the
morphism $\phi$ factors set-theoretically through the affinoid domain 
$\bS(a^d)$.  Hence $\phi$ induces a morphism
$\bS(a)\to\bS(a^d)$, so the homomorphism $K[s]\to K\angles{at\inv, t}$
extends to a homomorphism 
\[ F: K\angles{a^d s\inv, s}\To K\angles{at\inv, t}\qquad
s\mapsto t^d(1+g(t)),\quad a^d s\inv\mapsto (at\inv)^d(1+g(t))\inv. \]
Since $|g|_{\sup} < 1$, the induced map on canonical reductions is
\[ \td F: k[\sigma_1, \sigma_2]/(\sigma_1\sigma_2-\bar a{}^d) 
\To k[\tau_1,\tau_2]/(\tau_1\tau_2-\bar a)\qquad
\sigma_i\mapsto\tau_i^d \]
where $\sigma_1$ (resp.\ $\sigma_2,\tau_1,\tau_2$) is the residue of
$a^d s\inv$ (resp.\ $s, at\inv, t$).  Now $F$ is finite because
$\td F$ is finite~\cite[Theorem~6.3.5/1]{bgr:nonarch}, and 
it is easy to see that $F$ has degree $d$.  Flatness of $F$ is automatic
because its source 
and target are principal ideal domains: any affinoid algebra is noetherian,
and if $\sM(\cA)$ is an affinoid subdomain of 
$\A^{1,\berk} = \Spec(K[t])^\berk$ then any
maximal ideal of $\cA$ is the extension of a maximal ideal of $K[t]$
by~\cite[Lemma~5.1.2(1)]{conrad:irredcomps}.

For~(3), as above $\phi$ factors through $\bS(1,1)$ if we assume
$\alpha=1$, so we get a homomorphism
$F: K\angles{a^d s\inv, s}\to K\angles{t,t\inv}$.  In this case
the map $\td F$ on canonical reductions is clearly not finite, so $F$ is
not finite.\qed

\paragraph[The skeleton of a standard generalized annulus]
Define a section $\sigma:\R\to\G_m^\an$
of the tropicalization map $\trop:\G_m^\an\to\R$ by
\begin{equation} \label{eq:gauss.point.r}
\sigma(r) = \|\cdot\|_r \sptxt{where}
\left\|\sum_{n=-\infty}^\infty a_n t^n\right\|_r
= \max \big\{ |a_n|\cdot\exp(-rn) ~:~ n\in\Z \big\}. 
\end{equation}
When $r\in G$ the point $\sigma(r)$ is the Shilov boundary point of the
(strictly) affinoid domain $\trop\inv(r)$, and when $r\notin G$ we have
$\trop\inv(r) = \{\sigma(r)\}$.
The map $\sigma$ is easily seen to be continuous, and is in fact the only
continuous section of $\trop$.
We restrict $\sigma$ to obtain continuous sections
\[\begin{split} [\val(b),\val(a)] \To \bS(a,b) &\qquad
(\val(b),\val(a)) \To \bS(a,b)_+ \\
(\val(a),\infty) \To \bS(0,a)_+  &\qquad
(-\infty,\val(a)) \To \bS(a,\infty)_+
\end{split}\]
of $\trop$.

\begin{defn*}
  Let $A$ be a standard generalized annulus.  The
  \emph{skeleton} of $A$ is the closed subset
  \[ \Sigma(A)\coloneq\sigma(\R)\cap A = \sigma(\trop(A)). \]
  More explicitly, the skeleton of $\bS(a,b)$ (resp.\ $\bS(a,b)_+$, 
  resp.\ $\bS(0,a)_+$, resp.\ $\bS(a,\infty)_+$) is 
  \[\begin{split}
    \Sigma(\bS(a,b)) &\coloneq \sigma(\R)\cap\bS(a,b) 
    = \sigma([\val(b),\val(a)]) \\
    \Sigma(\bS(a,b)_+) &\coloneq \sigma(\R)\cap\bS(a,b)_+ 
    = \sigma((\val(b),\val(a))) \\
    \Sigma(\bS(0,a)_+) &\coloneq \sigma(\R)\cap\bS(0,a)_+
    = \sigma((\val(a),\infty))) \\
    \Sigma(\bS(a,\infty)_+) &\coloneq \sigma(\R)\cap\bS(a,\infty)_+
    = \sigma((-\infty,\val(a))).
  \end{split}\]
  We identify $\Sigma(A)$ with the interval/ray $\trop(A)$ via $\trop$ or
  $\sigma$. 
\end{defn*}

Note that $\tau_A\coloneq\sigma\circ\trop$ is a retraction of a standard generalized
annulus $A$ onto its skeleton.  This can be shown to be a strong deformation
retraction~\cite[Proposition~4.1.6]{berkovich:analytic_geometry}. 
Note also that the length of the skeleton of a standard generalized
annulus is equal to its modulus.

The set-theoretic skeleton has the following intrinsic characterization: 

\begin{prop}[{\cite[Proposition~2.2.5]{thuillier:thesis}}]
  \label{prop:skeleton.as.set}
  The skeleton of a standard generalized annulus is the set of all points that
  do not admit an affinoid neighborhood isomorphic to $\B(1)$. 
\end{prop}

\smallskip
The skeleton behaves well with respect to maps between standard generalized annuli:

\begin{prop} \label{prop:annulus.skeleton}
  Let $A$ be a standard generalized annulus of nonzero modulus 
  and let $\phi: A\to\G_m^\an$ be a morphism.  Suppose that 
  $\trop\circ\phi:\Sigma(A)\to\R$ is not constant.  Then:
  \begin{enumerate}
  \item For $x\in\Sigma(A)$ we have
    \[ \trop\circ\phi(x) = d\trop(x) + \val(\alpha) \]
    for some nonzero integer $d$ and some $\alpha\in K^\times$.  
  \item Let $B = \phi(A)$.  Then
    $B = \trop\inv(\trop(\phi(A)))$ is a standard generalized annulus in
    $\G_m^\an$ of the 
    same type, and $\phi:A\to B$ is a finite morphism of degree $|d|$.
  \item $\phi(\Sigma(A)) = \Sigma(B)$ and the following square commutes:
    \[ \xymatrix{
      {\trop(A)} \ar[rr]^{d(\cdot)+\val(\alpha)} \ar[d]_\sigma & &
      {\trop(B)} \ar[d]^\sigma \\
      {\Sigma(A)} \ar[rr]^\phi & &
      {\Sigma(B)} 
    }\]
  \end{enumerate}
\end{prop}

\pf Let $A'\cong\bS(a)\subset A$ be a standard closed annulus of nonzero
modulus such that $\trop\circ\phi$ is not constant on $\Sigma(A')$.  The
morphism $\phi$ is determined by a unit 
$f\in K\angles{at\inv, t}^\times$, and for $x\in\Sigma(A')$ we have
$\trop(\phi(x)) = -\log|f(x)|$.  Writing
$f(t) = \alpha\, t^d(1+g(t))$ as in~\eqref{eq:annulus.units}, if 
$r = \trop(x)$ then $-\log |f(x)| = -\log\|f\|_r = dr + \val(\alpha)$
since $\|1+g\|_r = 1$.  Since $\trop\circ\phi$ is nonconstant on 
$\Sigma(A')$ we must have $d\neq 0$.  Part~(1) follows by writing $A$ as
an increasing union of standard closed annuli and applying the same
argument.  The equality $B = \trop\inv(\trop(\phi(A)))$
follows from Proposition~\ref{prop:annulus.morphism}(2) in
the same way; since $\trop(\phi(A))$ is a closed interval (resp.\ open
interval, resp.\ open ray) when $\trop(A)$ is a closed interval (resp.\ open
interval, resp.\ open ray), it follows that $B$ is a standard generalized
annulus of the same type as $A$.

For part~(3) it suffices to show that
$\phi(\sigma(r)) = \sigma(dr+\val(\alpha))$ for 
$r\in\trop(A)$.  This follows from the above because
$\sigma(dr+\val(\alpha))$ is the supremum norm on
$\trop\inv(dr+\val(\alpha))$ (when $r\in G$) and 
$\phi$ maps $\trop\inv(r)$ surjectively onto 
$\trop\inv(dr+\val(\alpha))$.\qed

\begin{cor} \label{cor:finite.annulus.map}
  Let $\phi: A_1\to A_2$ be a finite morphism of standard generalized annuli and
  let $d$ be the degree of $\phi$.  
  Then $\phi(\Sigma(A_1)) = \Sigma(A_2)$, 
  $\phi(\sigma(r)) = \sigma(\pm dr+\val(\alpha))$
  for all $r\in\trop(A_1)$ and some $\alpha\in K^\times$, 
  and the modulus of $A_2$ is $d$ times the modulus of $A_1$.
  In particular, two standard generalized annuli of the same type are
  isomorphic if and only if they have the same modulus.
\end{cor}

\pf If the modulus of $A_1$ is zero then the result follows easily
from Proposition~\ref{prop:annulus.morphism}. Suppose that the modulus of
$A_1$ is nonzero.
By Proposition~\ref{prop:annulus.skeleton}, the only thing to show is that
$\trop\circ\phi$ is not constant on $\Sigma(A)$.  This is an immediate
consequence of Proposition~\ref{prop:annulus.morphism}(3).\qed

\paragraph[General annuli and balls]
In order to distinguish the properties of a standard generalized annulus
and its skeleton that are invariant under isomorphism, it is convenient to
make the following definition.

\begin{defn*}
  A \emph{closed ball} (resp.\ \emph{closed annulus}, 
  resp.\ \emph{open ball}, resp.\ \emph{open annulus}, 
  resp.\ \emph{punctured open ball}) is a $K$-analytic space isomorphic to 
  a standard closed ball (resp.\ standard closed annulus,
  resp. standard open ball, resp.\ standard open annulus, 
  resp.\ standard punctured open ball).  
  A \emph{generalized annulus} is a closed annulus, an open
  annulus, or a punctured open ball, and a \emph{generalized open annulus}
  is an open annulus or a punctured open ball.
\end{defn*}

\paragraph \label{par:generalized.annulus}
Let $A$ be a generalized annulus and fix an isomorphism
$\phi: A\isom A'$ with a standard generalized annulus $A'$.
The \emph{skeleton} of $A$ is defined to be 
$\Sigma(A) \coloneq \phi\inv(\Sigma(A'))$.
By Proposition~\ref{prop:skeleton.as.set} (or Corollary~\ref{cor:finite.annulus.map}) this
is a well-defined closed subset 
of $\Sigma(A)$.  We will view $\Sigma(A)$ as a closed interval
(resp.\ open interval, resp.\ open ray) with endpoints in $G$,
well-defined up to affine transformations of the form
$r\mapsto\pm r+\val(\alpha)$ for $\alpha\in K^\times$.  In particular
$\Sigma(A)$ is naturally a metric space, and it makes sense to talk about
piecewise affine-linear functions on $\Sigma(A)$ and of the slope of an affine-linear
function on $\Sigma(A)$ up to sign.  We remark that
$G\cap\Sigma(A)$ is equal to the set of type-$2$ points of $A$ contained
in $\Sigma(A)$.

The retraction $\tau_{A'} = \sigma\circ\trop: A'\to\Sigma(A')$ induces a
retraction $\tau_A: A\to\Sigma(A)$.  By Proposition~\ref{prop:annulus.skeleton} this
retraction is also independent of the choice of $A'$.

\begin{defn}
  Let $A$ be a generalized annulus, an open ball, or a closed ball.  A
  \emph{meromorphic function on $A$} 
  is by definition a quotient of an analytic function on $A$ by a nonzero
  analytic function on $A$.
\end{defn}

\smallskip
Note that a meromorphic function $f$ on $A$ is analytic
on the open analytic domain of $A$ obtained by deleting the poles
of $f$. If $A$ is affinoid then $f$ has only finitely many poles. 

Let $A$ be a generalized annulus, let $F:\Sigma(A)\to\R$ be a piecewise
affine function, and let $x$ be contained in the interior of $\Sigma(A)$.  
The \emph{change of slope of $F$ at $x$} is defined to be
\[ \lim_{\epsilon\to 0} \big( F'(x+\epsilon) - F'(x-\epsilon) \big); \]
this is independent of the choice of identification of $\Sigma(A)$ with
an interval in $\R$.  

We will need the following special case of the
Slope Formula~\parref{thm:PL}.
Its proof is an easy Newton polygon computation. 

\begin{prop} \label{prop:simple.PL}
  Let $A$ be a generalized annulus, let $f$ be a meromorphic function on
  $A$, and define $F:\Sigma(A)\to\R$ by $F(x) = -\log |f(x)|$.  
  \begin{enumerate}
  \item $F$ is a piecewise affine function with integer slopes, and
    for $x$ in the interior of $\Sigma(A)$ the change of slope of
    $F$ at $x$ is equal to the number of poles of $f$ retracting to $x$
    minus the number of zeros of $f$ retracting to $x$, counted with
    multiplicity. 

  \item Suppose that $A = \bS(0)_+$ and that $f$ extends to a meromorphic
    function on $\B(1)_+$.  Then for all $r\in(0,\infty)$ such that
    $r > \val(y)$ for all zeros and poles $y$ of $f$ in $A$, we have
    $F'(r) = \ord_0(f)$. 
  \end{enumerate}
\end{prop}

\begin{cor} \label{cor:decreasing.logf}
  Let $f$ be an analytic function on $\bS(0)_+$ that extends to a
  meromorphic function on $\B(1)_+$ with a pole at $0$ of order $d$.
  Suppose that $f$ has fewer than $d$ zeros on $\bS(0)_+$.
  Then $F = \log |f|$ is a monotonically increasing function on
  $\Sigma(\bS(0)_+) = (0,\infty)$.
\end{cor}

\smallskip
The following facts will also be useful:

\begin{lem} \label{lem:subtract.annulus.skeleton}
  Let $A$ be a generalized annulus.  Then the open analytic domain
  $A\setminus \Sigma(A)$ is isomorphic to an infinite disjoint union of
  open balls.  Each connected component $B$ of $A\setminus\Sigma(A)$
  retracts onto a single point $x\in\Sigma(A)$, and the closure of $B$ in $A$ is
  equal to $B\cup\{x\}$.
\end{lem}

\pf First we assume that $A$ is the standard closed annulus 
$\bS(1)=\sM(K\angles{t^{\pm 1}})$ of modulus
zero.  Then $\Sigma(A) = \{x\}$ is the Shilov boundary point of $A$.  The
canonical reduction of $A$ is isomorphic to $\G_{m,k}$, the inverse image
of the generic point of $\G_{m,k}$ is $x$, the
inverse image of a residue class $\bar y\in k^\times = \G_m(k)$ is the open ball
$\{\|\cdot\|~:~\|t-y\|<1\}$ (where $y\in R^\times$ reduces to $\bar y$),
and the fibers over the closed points of 
$\G_{m,k}$ are the connected components of $A\setminus\{x\}$
by~\cite[Lemme~2.1.13]{thuillier:thesis}.  This proves the first
assertion, and the second follows from the anti-continuity of the
reduction map.

Now let $A$ be any generalized annulus; we may assume that 
$A$ is standard.  Let $r\in \trop(A)$.  If $r\notin G$ then
$\trop\inv(r)$ is a single point of type $3$, so suppose
$r\in G$, say $r = \val(a)$ for $a\in K^\times$.  After translating by
$a\inv$ we may and do assume that $r = 0$, so 
$\trop\inv(r) = \bS(1)= \sM(K\angles{t^{\pm 1}})$. 
The subset $\bS(1)\setminus\{\sigma(0)\}$ is clearly closed in
$A\setminus\Sigma(A)$, and it is open as well since it is the union of
the open balls $\{\|\cdot\|~:~\|t-y\|<1\}$ for $y\in R^\times$.  Therefore
the connected components of $\bS(1)\setminus\{\sigma(0)\}$ are
also connected components of $A\setminus\Sigma(A)$, so we are reduced to
the case treated above.\qed 

\begin{lem} \label{lem:unit.logf.constant}
  Let $A$ be a generalized annulus and let $f$ be a unit on $A$.  Then
  $x\mapsto\log|f(x)|$ factors through the retraction
  $\tau_A:A\to\Sigma(A)$.  In particular, $x\mapsto\log|f(x)|$ is locally
  constant away from $\Sigma(A)$.
\end{lem}

\pf This follows immediately from Lemma~\ref{lem:subtract.annulus.skeleton}
and the elementary fact that a unit on an open ball has constant absolute
value.\qed

\section{Semistable decompositions and skeleta of curves}
\label{sec:skeleta.curves}

For the rest of this paper $X$ denotes a smooth connected algebraic curve over $K$,
$\hat X$ denotes its smooth completion, and $D = \hat X\setminus X$
denotes the set of punctures.  We will define a skeleton
inside of $X$ relative to the following data.

\begin{defn}
  A \emph{semistable vertex set of $\hat X$} is a finite set $V$ of type-$2$
  points of $\hat X^\an$ such that $\hat X^\an\setminus V$ is a disjoint union of
  open balls and finitely many open annuli.
  A \emph{semistable vertex set of $X$} is a semistable vertex set of
  $\hat X$ such that the punctures in $D$ are contained in distinct
  connected components of $\hat X^\an\setminus V$ isomorphic to open balls.
  A decomposition of $X^\an$ into a semistable vertex set and a disjoint
  union of open balls and finitely many generalized open annuli is called a
  \emph{semistable decomposition of $X$}.  
\end{defn}

When we refer to `an open ball in a semistable decomposition of $X$' or
`a generalized open annulus in a semistable decomposition of $X$' we
will always mean a connected component of $X^\an\setminus V$ of the
specified type.
Note that the punctured open balls in a semistable
decomposition of $X$ are in bijection with $D$, and that
there are no punctured open balls in a
semistable decomposition of a complete curve.  A semistable
vertex set of $X$ is also a semistable vertex set of $\hat X$.  

The semistable vertex sets of $\hat X$ correspond naturally and
bijectively to isomorphism classes of semistable formal models of $\hat
X$.  See Theorem~\ref{thm:construct.models}.

\begin{lem} \label{lem:skeleton.is.graph}
  Let $V$ be a semistable vertex set of $X$, let $A$ be a connected
  component of $X^\an\setminus V$, and let $\bar A$ be the closure of $A$
  in $\hat X^\an$.  Let $\del_{\lim} A = \bar A\setminus A$ be the
  \emph{limit boundary} of $A$, i.e.\ the set of limit points of $A$ in
  $\hat X^\an$ that are not contained in $A$.%
  \footnote{As opposed to the canonical boundary discussed in 
    \cite[\S{2.5.7}]{berkovich:analytic_geometry}.}
  \begin{enumerate}
  \item If $A$ is an open ball then $\del_{\lim} A = \{x\}$ for some 
    $x\in V$.
  \item Suppose that $A$ is an open annulus, and fix an isomorphism
    $A\cong\bS(a)_+$.  Let $r=\val(a)$.  Then $\sigma:(0,r)\to A$ extends
    in a unique way to a continuous map $\sigma:[0,r]\to X^\an$ such that
    $\sigma(0),\sigma(r)\in V$, and 
    $\del_{\lim}A = \{\sigma(0),\sigma(r)\}$.  (It may happen that 
    $\sigma(0) = \sigma(r)$.)
  \item Suppose that $A$ is a punctured open ball, and fix an isomorphism 
    $A\cong\bS(0)_+$.  Then $\sigma:(0,\infty)\to A$ extends in a unique
    way to a continuous map $\sigma:[0,\infty]\to\hat X^\an$ such that
    $\sigma(0)\in V$, $\sigma(\infty)\in D$,
    and $\del_{\lim} A = \{\sigma(0),\sigma(\infty)\}$.
  \end{enumerate}
\end{lem}

\pf First note that in~(1) and~(2), $\bar A$ is the closure of $A$ in 
$X^\an$ because every point of $\hat X\setminus X$ has an open
neighborhood disjoint from $A$.  Since $A$ is closed in $X\setminus V$,
its limit boundary is contained in $V$.  

Suppose that $A$ is an open ball, and fix an isomorphism 
$\phi:\B(1)_+\isom A$.
For $r\in(0,\infty)$ we define $\|\cdot\|_r\in\B(1)_+$
by~\eqref{eq:gauss.point.r}.  Fix an affine open subset $X'$ of $X$ such
that $A\subset(X')^\an$.  For any $f\in K[X']$
the map $r\mapsto\log\|f\|_r$ is piecewise affine with
finitely many changes in slope by Proposition~\ref{prop:simple.PL}.  Therefore we
may define $\|f\| = \lim_{r\to 0} \|f\|_r\in\R$.  The map
$f\mapsto\|f\|$ is easily seen to be a multiplicative norm on $K[X']$,
hence defines a point $x\in(X')^\an\subset X^\an$.

Let $y$ be the Shilov point of $\B(1)$ and let $A' = \B(1)_+\cup\{y\}$.
Since $\B(1)\setminus\{y\}$ is a disjoint union of open balls it is clear
that $A'$ is a closed, hence compact subset of $\B(1)$.  Extend $\phi$ to
a map $A'\to(X')^\an\subset X^\an$ by $\phi(y) = x$.  We claim that $\phi$ is
continuous.  By the definition of the topology on $(X')^\an$ it suffices
to show that the set $U = \{ z\in A'~:~|f(\phi(z))|\in(c_1,c_2)\}$ is open
for all $f\in K[X']$ and all $c_1<c_2$.  Since $U\cap\B(1)_+$ is open, we
need to show that $U$ contains a neighborhood of $y$ if $y\in U$, 
i.e.,\ if $\|f\|\in(c_1,c_2)$.  Choose $a\in\fm_R\setminus\{0\}$ such that
$f$ has no zeros in $\bS(a)_+$.  Note that $\bS(a)_+\cup\{y\}$ is a
neighborhood of $y$ in $A'$.  Since $\|f\| = \lim_{r\to 0} \|f\|_r$ we
have that $\|f\|_r\in(c_1,c_2)$ for $r$ close enough to $1$; hence we may
shrink $\bS(a)_+$ so that $\phi(\Sigma(\bS(a)_+))\subset(r_1,r_2)$.
With Lemma~\ref{lem:unit.logf.constant} this implies that
$\bS(a)_+\subset U$, so $\phi$ is indeed continuous.  Since
$A'$ is compact we have that $\phi(A') = A\cup\{x\}$ is closed, which
completes the proof of~(1).

If $A\cong\bS(0)_+$ is a punctured open ball then certainly the puncture
$0$ is in $\bar A$.  The above argument effectively proves the rest
of~(3), and~(2) is proved in exactly the same way.\qed

\begin{defn}
  Let $V$ be a semistable vertex set of $X$.
  The \emph{skeleton} of $X$ with respect to $V$ is
  \[ \Sigma(X,V) = V\cup\bigcup\Sigma(A) \]
  where $A$ runs over all of the connected components of 
  $X^\an\setminus V$ that are generalized open annuli.
\end{defn}

\begin{lem} \label{lem:skeleton.nice}
  Let $V$ be a semistable vertex set of $X$ and let $\Sigma=\Sigma(X,V)$ be
  the associated skeleton.  Then:
  \begin{enumerate}
  \item $\Sigma$ is a closed subset of $X^\an$ which is compact if and only if 
    $X = \hat X$.
  \item The limit boundary of $\Sigma$ in $\hat X^\an$ is equal to
    $D$. 
  \item The connected components of $X^\an\setminus\Sigma(X,V)$ are open
    balls, and the limit boundary $\del_{\lim} B$ of any connected
    component $B$ is a single point $x\in\Sigma(X,V)$.
  \item $\Sigma$ is equal to the set of points in $X^\an$ that
    do not admit an affinoid neighborhood isomorphic to $\B(1)$ and
    disjoint from $V$.
  \end{enumerate}
\end{lem}

\pf The first two assertions are clear from Lemma~\ref{lem:skeleton.is.graph},
and the third follows from Lemmas~\ref{lem:skeleton.is.graph}
and~\ref{lem:subtract.annulus.skeleton}.
Let $\Sigma'$ be the set of points in $X^\an$ that do not admit an
affinoid neighborhood isomorphic to $\B(1)$ and disjoint from $V$.
We have $\Sigma'\subset\Sigma$ by~(3).
For the other inclusion, let $x\in\Sigma$.  If $x\in V$ then clearly
$x\in\Sigma'$, so suppose $x\notin V$.  Then the connected component $A$ of
$x$ in $X^\an\setminus V$ is a generalized open annulus; since any
connected neighborhood of $x$ is contained in $A$, we have $x\in\Sigma'$
by Proposition~\ref{prop:skeleton.as.set}.\qed

\begin{defn} \label{defn:completed.skeleton}
  Let $V$ be a semistable vertex set of $X$.  The \emph{completed skeleton}
  of $X$ with respect to $V$ is defined to be the closure of 
  $\Sigma(X,V)$ in $\hat X^\an$ and is denoted $\hat\Sigma(X,V)$, so
  $\hat\Sigma(X,V) = \Sigma(X,V)\cup D$.  The completed skeleton has the
  structure of a graph with vertices $V\cup D$; the interiors of the edges
  of $\Sigma(X,V)$ are the skeleta of the generalized open annuli in the
  semistable decomposition of $X$ coming from $V$.
\end{defn}

\begin{rem} \label{rem:skel.is.graph}
  By Lemma~\ref{lem:skeleton.nice}(1), if $X=\hat X$ then the skeleton
  $\Sigma(X,V) = \hat\Sigma(X,V)$ is a finite metric graph
  (cf.~\parref{par:dual.graph}).  If $X$ is not proper then 
  $\hat\Sigma(X,V)$ is a finite graph with vertex set $V\cup D$,
  but the ``metric'' on $\hat\Sigma(X,V)$ is degenerate since it has edges
  of infinite length.   See~\parref{defn:abstract.complex}. 
\end{rem}

\begin{defn}
  Let $V$ be a semistable vertex set of $X$ and let 
  $\Sigma = \Sigma(X,V)$.  We define a retraction
  $\tau_V = \tau_\Sigma :X^\an\to\Sigma$ as follows.  Let 
  $x\in X^\an\setminus\Sigma$ and let $B_x$ be the connected component
  of $x$ in $X^\an\setminus\Sigma$.  Then $\del_{\lim}(B_x) = \{y\}$
  for a single point $y\in X^\an$; we set $\tau_V(x) = y$.
\end{defn}

\begin{lem} \label{lem:tau.continuous}
  Let $V$ be a semistable vertex set of $X$.  The retraction
  $\tau_V: X^\an\to\Sigma(X,V)$ is continuous, and if $A$ is 
  a generalized open annulus in the semistable decomposition of
  $X$ then $\tau_V$ restricts to the retraction
  $\tau_A: A\to\Sigma(A)$ defined in~\parref{par:generalized.annulus}. 
\end{lem}

\pf The second assertion follows
from Lemma~\ref{lem:subtract.annulus.skeleton}, so $\tau_V$ is continuous when
restricted to any connected component $A$ of $X^\an\setminus V$ which is
a generalized open annulus.  Hence it is enough to show that if 
$x\in V$ and $U$ is an open neighborhood of $x$ then $\tau_V\inv(U)$
contains an open neighborhood of $x$.  This is left as an exercise to the
reader.\qed 

\begin{prop}
  Let $V$ be a semistable vertex set of $X$.  Then 
  $\Sigma(X,V)$ and $\hat\Sigma(X,V)$ are connected.
\end{prop}

\pf This follows from the continuity of $\tau_V$
and the connectedness of $X^\an$.\qed

\smallskip
The skeleton of a curve naturally carries the following kind of
combinatorial structure, which is similar to that of a metric graph.

\begin{defn} \label{defn:abstract.complex}
  A \emph{dimension-$1$ abstract $G$-affine polyhedral complex} is a
  combinatorial object $\Sigma$ consisting of the following data.  We are
  given a finite discrete set $V$ of \emph{vertices} and a collection of
  finitely many \emph{segments} and \emph{rays}, where a segment is 
  a closed interval in $\R$ with distinct endpoints in $G$ and a ray is a
  closed ray in $\R$ with endpoint in $G$.  Segments and rays are only
  defined up to isometries of $\R$ of the form $r\mapsto\pm r + \alpha$ for
  $\alpha\in G$.  The segments and rays are collectively called \emph{edges} of
  $\Sigma$.  Finally, we are given an identification of the endpoints of the
  edges of $\Sigma$ with vertices.  The complex $\Sigma$ has
  an obvious realization as a topological space, which we will also denote
  by $\Sigma$.   If $\Sigma$ is connected then it is a metric space under
  the shortest-path metric.

  A \emph{morphism} of dimension-$1$ abstract $G$-affine polyhedral
  complexes is a continuous function $\phi:\Sigma\to\Sigma'$ sending
  vertices to vertices and such that if $e\subset\Sigma$ is an edge then
  either $\phi(e)$ is a vertex of $\Sigma'$, or 
  $\phi(e)$ is an edge of $\Sigma'$ and for all $r\in e$ we have
  $\phi(r) = dr+\alpha$ for a nonzero integer $d$ and some $\alpha\in G$.

  A \emph{refinement} of a dimension-$1$ abstract $G$-affine polyhedral
  complex is a complex $\Sigma'$ obtained from $\Sigma$ by inserting
  vertices at $G$-points of edges of $\Sigma$ and dividing those edges in
  the obvious way.  Note that $\Sigma$ and $\Sigma'$ have the same
  topological and metric space realizations.
\end{defn}

\begin{rem}
  Abstract integral $G$-affine polyhedral complexes of arbitrary dimension are defined
  in~\cite[\S~1]{thuillier:thesis} in terms of groups of integer-slope
  $G$-affine functions.  In the one-dimensional case
  the objects of loc.\ cit.\ are roughly the same as the
  dimension-$1$ abstract integral $G$-affine polyhedral complexes in the sense of
  our ad-hoc definition above, since the knowledge of what functions on a line
  segment have slope one is basically the same as the data of a metric.
  We choose to use this definition for concreteness and in order to
  emphasize the metric nature of these objects. 
\end{rem}

\paragraph \label{par:define.complex.structure}
Let $V$ be a semistable vertex set of $X$.  Then $\Sigma(X,V)$ is a
dimension-$1$ abstract $G$-affine polyhedral complex with vertex set $V$
whose edges are the closures of the skeleta of the generalized open annuli
in the semistable decomposition of $X$.  In particular, $\Sigma(X,V)$ is a
metric space, and each edge $e$ of $\Sigma(X,V)$ is identified via a local
isometry with the skeleton of the corresponding generalized open annulus.
Note that if $e$ is a segment then the length of $e$ is equal to the
modulus of the corresponding open annulus.  The $G$-points of 
$\Sigma(X,V)$ are exactly the type-$2$ points of $X$ contained in 
$\Sigma(X,V)$.

\begin{prop} \label{prop:relations.vertex.sets}
  Let $V$ be a semistable vertex set of $X$ and let $X'$ be a nonempty open
  subscheme of $X$.
  \begin{enumerate}
  \item Let $V'$ be a semistable vertex set of $X'$ containing $V$.  Then
    $\Sigma(X,V)\subset\Sigma(X',V')$ and $\Sigma(X',V')$ induces a
    refinement of $\Sigma(X,V)$.  Furthermore,
    $\tau_{\Sigma(X,V)}\circ\tau_{\Sigma(X',V')} = \tau_{\Sigma(X,V)}$.
  \item Let $V'\subset\Sigma(X,V)$ be a finite set of type-$2$ points.
    Then $V\cup V'$ is a semistable vertex set of $X$ and 
    $\Sigma(X,V\cup V')$ is a refinement of $\Sigma(X,V)$.
  \item Let $W\subset X^\an$ be a finite set of type-$2$ points.  Then
    there is a semistable vertex set $V'$ of $X'$ containing 
    $V\cup W$.
  \end{enumerate}
\end{prop}

\pf In~(1), the inclusion $\Sigma(X,V)\subset\Sigma(X',V')$ follows
from Lemma~\ref{lem:skeleton.nice}(4), and the fact that $\Sigma(X',V')$
induces a refinement of $\Sigma(X,V)$ is an easy consequence of the
structure of morphisms of generalized open
annuli (Proposition~\ref{prop:annulus.skeleton}).  The equality
$\tau_{\Sigma(X,V)}\circ\tau_{\Sigma(X',V')} = \tau_{\Sigma(X,V)}$
follows from the definitions.  In~(2) we may assume that
$V'$ is a single point; it then suffices to show that if $A$ is a
generalized open annulus and $x\in\Sigma(A)$ is a type-$2$ point then
$A\setminus\{x\}$ is a disjoint union of generalized open annuli and open
balls.  Choose an identification of $A$ with a standard
generalized open annulus $\bS(a,b)_+$ such that 
$0\in(\trop(b),\trop(a))$ and $x = \sigma(0)$.  
As in the proof of Lemma~\ref{lem:subtract.annulus.skeleton} 
we have that $\bS(1)\setminus\{x\}$ is a disjoint union of open balls,
and is open and closed in $A$; hence 
\[ A\setminus\{x\} = \bS(a,1)_+\djunion (\bS(1)\setminus\{x\})\djunion
\bS(1,b)_+ \]
is a disjoint union of generalized open annuli and open balls.

It suffices to prove~(3) when $W = \{x\}$ and
$X=X'$, and when $W=\emptyset$ and $X\setminus X' = \{y\}$.
In the first case, we may assume that $x\notin\Sigma(X,V)$ by the above.
Suppose that the connected component $A$ of
$X^\an\setminus V$ containing $x$ is an open ball.  One shows as in~(2)
that $A\setminus\{x\}$ is a disjoint union of open balls and an open
annulus, so $V\cup\{x\}$ is a semistable vertex set.  If
$A$ is a generalized open annulus then the connected component of $x$ in
$A\setminus\{\tau_A(x)\}$ is an open ball, so $V\cup\{x,\tau_A(x)\}$ is a
semistable vertex set.  In the case $X\setminus X' = \{y\}$ one proceeds in
exactly the same way.\qed 

A semistable vertex set $V$ is called \emph{strongly semistable} if the
graph $\hat\Sigma(X,V)$ has no loop edges.  (See
Definition~\ref{defn:completed.skeleton}.)

\begin{cor}
  Any semistable vertex set of $X$ is contained in a strongly semistable
  vertex set of $X$.  
\end{cor}

\section{Relation with semistable models}

Recall that $X$ is a smooth connected algebraic curve over $K$,
$\hat X$ is its smooth completion, and $D = \hat X\setminus X$
is the set of punctures.  
The (formal) semistable reduction theory of a smooth complete algebraic curve was
worked out carefully in~\cite{bosch_lutk:uniformization} in the language
of rigid analytic spaces and formal analytic varieties (see
Remark~\ref{rem:remarks.semistable}(3)); one can view 
much of this section as a translation of that paper into our
language of semistable vertex sets. 

\paragraph
It is more natural in the context of analytic geometry to use 
Bosch-L\"utkebohmert's theory~\cite{bl:fmI} of admissible
formal models of $\hat X^\an$ instead of algebraic models of $\hat X$.
An \emph{admissible $R$-algebra} is an $R$-flat quotient of a
convergent power series ring 
\[ R\angles{X_1,\ldots,X_n} = \bigg\{\sum a_I X^I\in R\ps{X_1,\ldots,X_n} 
~:~ |a_I| \to 0 \text{ as } |I|\to\infty \bigg\} \]
by a finitely generated ideal.  An \emph{admissible formal $R$-scheme} is
a formal scheme $\fX$, locally modeled on the formal spectrum of an admissible 
$R$-algebra.  If $A$ is an admissible $R$-algebra then $A\tensor_R K$ is
a strictly $K$-affinoid algebra, and the functor taking $A$ to its
Berkovich spectrum $\sM(A\tensor_R K)$ glues to give the 
\emph{Raynaud generic fiber functor} $\fX\mapsto\fX^\an$ from the category
of admissible formal $R$-schemes to the category of $K$-analytic spaces.
There is a canonical anti-continuous reduction map
$\red: \fX^\an\to\fX_k$.

\begin{defn*}
  \begin{enumerate}
  \item A connected reduced algebraic curve over a field $k$ is called 
    \emph{semistable} if its singularities are ordinary double points.  It is
    \emph{strongly semistable} if in addition its irreducible components
    are smooth.
  \item A \emph{(strongly) semistable formal $R$-curve} is an integral
    admissible formal $R$-curve $\fX$ whose special fiber is
    a (strongly) semistable curve.  A 
    \emph{(strongly) semistable formal model} for
    $\hat X$ is a (strongly) semistable proper formal $R$-curve $\fX$
    equipped with an isomorphism $\fX^\an\cong\hat X^\an$.
  \end{enumerate}
\end{defn*}

\begin{rem}\label{rem:remarks.semistable}
  \begin{enumerate}
  \item Let $\fX$ be a (strongly) semistable formal $R$-curve.  Then
    $\fX$ is proper if and only if $\fX^\an$ is proper
    by~\cite[Corollary~4.4]{temkin:local_properties}.  Therefore the
    properness hypothesis in the definition of a (strongly) semistable
    formal model for $\hat X$ is redundant.

  \item A \emph{(strongly) semistable algebraic model} for $\hat X$ is a
    flat and integral proper relative curve $\cX\to\Spec(R)$ whose special
    fiber is a (strongly) semistable curve and whose generic fiber is
    equipped with an isomorphism to $\hat X$.  A (strongly) semistable
    algebraic model $\cX$ for $\hat X$ gives rise to a
    (strongly) semistable formal model $\fX$ for $\hat X$ by 
    completing.  Indeed, $\fX$ and $\cX$ have the same special fiber, and
    $\fX^\an\cong\hat X^\an$ by~\cite[Theorem~5.3.1(4)]{conrad:irredcomps}.
    Conversely, a (strongly) semistable formal model for
    $\hat X$ uniquely algebraizes to a (strongly) semistable
    algebraic model by a suitable formal GAGA theorem over
    $R$~\cite[Corollaire~2.3.19]{abbes:egrI}. 
    Hence there is no essential difference between the algebraic and
    formal semistable reduction theories of $\hat X$.  

  \item Let $\fX$ be a semistable formal $R$-curve.  Since $\bar\fX$ is
    reduced, $\fX$ is a formal analytic variety in the sense
    of~\cite{bosch_lutk:uniformization}.
    In particular, if $\Spf(A)$ is a formal affine open subset of $\fX$
    then $A$ is the ring of power-bounded elements of $A\tensor_R K$ ---
    that is, $\Spf(A)$ is the canonical model of $\sM(A\tensor_R K)$ ---
    and $A\tensor_R k$ is the canonical reduction of $A\tensor_R K$. 

  \end{enumerate}
\end{rem}

Let $a\in R\setminus\{0\}$.
The \emph{standard formal annulus of modulus $\val(a)$} is defined to be
\[ \fS(a) \coloneq \Spf\big(R\angles{s,t}/(st-a)\big). \]
This is the canonical model of the standard closed annulus $\bS(a)$.

\begin{prop} \label{prop:definitions.semistable}
  Let $\fX$ be a strongly semistable formal $R$-curve and let
  $\xi\in\fX$ be a singular point of $\bar\fX$.  There is a formal
  neighborhood $\fU$ of $\xi$ and an \'etale morphism
  $\phi:\fU\to\fS(a)$ for some $a\in\fm_R\setminus\{0\}$ such that 
  $\phi^\an$ restricts to an isomorphism $\red\inv(\xi)\isom\bS(a)_+$.
\end{prop}

\pf This is essentially~\cite[Proposition~2.3]{bosch_lutk:uniformization};
here we explain how the proof of loc.\ cit.\ implies the
proposition.  Shrinking $\fX$ if necessary, we may and do assume that 
$\fX = \Spf(A)$ is affine and connected, and that the maximal ideal
$\fm_\xi\subset\bar A = A\tensor_R k$ corresponding to $\xi$ is generated
by two functions  $\bar f,\bar g\in\bar A$ whose product is zero
(this is possible because $\fX$ is \emph{strongly} semistable).  Let
$k[x,y]/(xy)\to\bar A$ be the homomorphism  
sending $x\mapsto\bar f$ and $y\mapsto\bar g$.
Since $\xi$ is an ordinary double point, we can
choose $\bar f$ and $\bar g$ such that the map on completed local rings
$k\ps{x,y}/(xy)\to\hat\sO_{\bar\fX,\xi}$ is an isomorphism.  
It follows from~\cite[Prop. 17.6.3]{egaIV_4} that the morphism
$\Spec(\bar A)\to\Spec(k[x,y]/(xy))$ is \'etale at $\xi$, so shrinking $\fX$
further we may assume that $\Spec(\bar A)\to\Spec(k[x,y]/(xy))$ is \'etale.
One then proceeds as in the proof
of~\cite[Proposition~2.3]{bosch_lutk:uniformization} to find lifts
$f,g\in A$ of $\bar f,\bar g$ such that $fg = a\in R\setminus\{0\}$; the
induced morphism $\fX\to\fS(a)$ is \'etale because it lifts an \'etale
morphism on the special fiber.  The fact that $\phi$ restricts to an isomorphism
$\red\inv(\xi)\isom\fS(a)_+$ is part~(i) of loc.\ cit.\qed

\smallskip
The following characterization of strongly semistable formal $R$-curves is
also commonly used in the literature, for example
in~\cite[D\'efinition~2.2.8]{thuillier:thesis} (see also
Remarque~2.2.9 in loc.\ cit.).  

\begin{cor} 
  An integral admissible formal $R$-curve $\fX$ is strongly
  semistable if and only if it has a covering by Zariski-open sets $\fU$
  which admit an \'etale morphism to $\fS(a_\fU)$ for some 
  $a_\fU\in R\setminus\{0\}$.
\end{cor}

\smallskip
Let $\fX$ and $\fX'$ be two semistable formal models for $\hat X$.  We say that 
$\fX$ \emph{dominates} $\fX'$, and we write $\fX\geq\fX'$, if there exists
an $R$-morphism $\fX\to\fX'$ inducing the identity on the generic fiber 
$\hat X^\an$.  Such a morphism is unique if it exists.  The relation $\geq$ is
a partial ordering on the set of semistable formal models for $\hat X$.  
(We will always consider semistable formal models of $\hat X$ up to
isomorphism; any isomorphism is unique.)

\paragraphtoc[Semistable models and semistable decompositions]
The special fiber of the canonical model for $\B(1)$ is isomorphic to 
$\A^1_k$, and the inverse image of the origin is the open unit ball
$\B(1)_+$.  When $|a| < 1$ the special fiber of $\fS(a)$ is isomorphic to
$k[x,y]/(xy)$, and the inverse image of the
origin under the reduction map is $\bS(a)_+$.  
The following much stronger version of these facts provides the relation
between semistable models and semistable decompositions of $\hat X$.

\begin{thm}[Berkovich, Bosch-L\"utkebohmert] \label{thm:inverse.image.reduction}
  Let $\fX$ be an integral admissible  formal $R$-curve with
  reduced special fiber and let $\xi\in\bar\fX$ be any point.
  \begin{enumerate}
  \item $\xi$ is a generic point if and only if
    $\red\inv(\xi)$ is a single type-$2$ point of $\fX^\an$.
  \item $\xi$ is a smooth closed point if and only if
    $\red\inv(\xi)\cong\B(1)_+$.
  \item $\xi$ is an ordinary double point if and only if
    $\red\inv(\xi)\cong\bS(a)_+$ for some $a\in\fm_R\setminus\{0\}$.
  \end{enumerate}
\end{thm}

\pf As in Remark~\ref{rem:remarks.semistable}(3) the hypothesis on the
special fiber of $\fX$ allows us to view $\fX$ as a formal analytic
variety.  Hence the first statement follows
from~\cite[Proposition~2.4.4]{berkovich:analytic_geometry}, and
the remaining assertions 
are~\cite[Propositions~2.2 and 2.3]{bosch_lutk:uniformization}.\qed

\smallskip
Let $\fX$ be a semistable formal model for $\hat X$.  We let
$V(\fX)$ denote the inverse image of the set of generic points of $\bar\fX$
under the reduction map.
This is a finite set of type-$2$ points of $\fX^\an$ that maps bijectively
onto the set of generic points of $\bar\fX$.

\begin{cor}
  Let $\fX$ be a semistable formal model for $\hat X$.
  Then $V(\fX)$ is a semistable vertex set of $\hat X$, and the decomposition of 
  $\hat X^\an\setminus V(\fX)$ into formal fibers is a semistable
  decomposition. 
\end{cor}

\pf By~\cite[Lemme~2.1.13]{thuillier:thesis} the formal fibers of $\fX$
are the connected components of $\hat X^\an\setminus V(\fX)$, so the assertion
reduces to Theorem~\ref{thm:inverse.image.reduction}.\qed

\begin{rem}
  The semistable vertex set $V(\fX)$ is a semistable vertex set of the
  punctured curve $X = \hat X\setminus D$ if and only
  if the punctures $x\in D$ reduce to distinct smooth closed points of
  $\fX(k)$, that is, if and only if $\fX$ is a semistable model of the
  marked curve $(\hat X,D)$.
\end{rem}

\paragraph \label{par:dual.graph}
Let $\fX$ be a semistable formal model for $\hat X$.  Let $\xi\in\bar\fX$
be a singular point and let $z_1,z_2\in\hat X^\an$ be the inverse images
of the generic points of $\bar\fX$ specializing to $\xi$ (it may be that $z_1=z_2$).
Then $z_1,z_2$ are the vertices of the edge in $\Sigma(\hat X,V(\fX))$
whose interior is $\Sigma(\red\inv(\xi))$ by the anti-continuity of the
reduction map and Lemma~\ref{lem:skeleton.is.graph}(2).  It follows that 
$\Sigma(\hat X,V(\fX))$ is the incidence graph of $\bar\fX$
(cf.\ Remark~\ref{rem:skel.is.graph}).
In other words, the vertices of  
$\Sigma(\hat X,V(\fX))$ correspond to irreducible components of $\bar\fX$ and the
edges of $\Sigma(\hat X,V(\fX))$ correspond to the points where the components of
$\bar\fX$ intersect.  Moreover, if $\fX$ admits an \'etale map to some
$\fS(a) = \Spf(R\angles{x,y}/(xy-a))$ in a neighborhood of $\xi$, then
$\val(a)$ is the length of the edge corresponding to $\xi$ (see the proof
of Proposition~\ref{prop:skeleta.agree}). 

It is clear from the above that a semistable formal model
$\fX$ for $\hat X$ is strongly semistable if and only if $V(\fX)$ is a
strongly semistable vertex set.

Berkovich~\cite{berkovich:locallycontractible2} and
Thuillier~\cite{thuillier:thesis} define the skeleton 
of a \emph{strongly} semistable formal $R$-curve
using Proposition~\ref{prop:definitions.semistable}.  In order to use their
results, we must show that the two notions of the skeleton agree:

\begin{prop} \label{prop:skeleta.agree}
  Let $\fX$ be a strongly semistable formal model for $\hat X$.  The skeleton
  $\Sigma(\hat X, V(\fX))$ is naturally identified with the skeleton of 
  $\fX$ defined in~\cite{thuillier:thesis} as dimension-$1$ abstract
  $G$-affine polyhedral complexes.
\end{prop}

\pf Thuillier~\cite[D\'efinition~2.2.13]{thuillier:thesis} defines the
skeleton $S(\fX)$ of $\fX$
to be the set of all points that do not admit an affinoid neighborhood
isomorphic to $\B(1)$ and disjoint from $V(\fX)$, so $\Sigma(\hat X, V(\fX))=S(\fX)$
as sets by Lemma~\ref{lem:skeleton.nice}(4).  Let
$\xi\in\bar\fX$ be a singular point and let $\fU$ be a formal affine
neighborhood of $\xi$ admitting an \'etale morphism $\phi:\fU\to\fS(a)$ and
inducing an isomorphism $\red\inv(\xi)\isom\bS(a)_+$ as
in Proposition~\ref{prop:definitions.semistable}.  Shrinking $\fU$ if necessary, we
may and do assume that $\xi$ is the only singular point of $\bar\fU$ and
that $\bar\fU$ has two generic points $\zeta_1,\zeta_2$.  Let
$z_1,z_2\in V(\fX)$ be the inverse images of $\zeta_1,\zeta_2$.  Then
$\Sigma(X,V(\fX))\cap\fU^\an$ is the edge in $\Sigma(X,V(\fX))$ connecting
$z_1,z_2$ with interior $\Sigma(\red\inv(\xi))$.  Since
$\phi^\an$ maps $\red\inv(\xi)$ isomorphically onto $\bS(a)_+$ it induces an
isometry $\Sigma(X,V(\fX))\cap\fU^\an\isom\Sigma(\bS(a))$.  
The polyhedral structure on $S(\fX)\cap\fU^\an$ is more or less by
definition induced by the identification of $\Sigma(\bS(a))$ with
$[0,\val(a)]$; see~\cite[Th\'eor\`eme~2.2.10]{thuillier:thesis}.  Hence 
$\Sigma(\hat X,V(\fX))=S(\fX)$ as $G$-affine polyhedral complexes.\qed

\smallskip
In order to prove that semistable vertex sets are in one-to-one
correspondence with semistable models as above, it remains to construct a
semistable model from a semistable decomposition.  The following theorem
is folklore; while it is well-known to experts, and in
some sense is implicit in~\cite{temkin:stable_modification}, we have been
unable to find an explicit reference.

\begin{thm} \label{thm:construct.models}
  The association $\fX\mapsto V(\fX)$ sets up a bijection between the set
  of semistable formal models of $\hat X$ (up to isomorphism) and the set
  of semistable vertex   sets of $\hat X$.  Furthermore, 
  $\fX$ dominates $\fX'$ if and only if $V(\fX')\subset V(\fX)$.
\end{thm}

\smallskip
We will need the following lemmas in the proof
of Theorem~\ref{thm:construct.models}. 

\begin{lem} \label{lem:delete.open.affinoid}
  \begin{enumerate}
  \item Let $B\subset\hat X^\an$ be an analytic open subset isomorphic to 
    an open ball.  Then $\hat X^\an\setminus B$ is an affinoid domain in
    $\hat X^\an$. 
  \item Let $A\subset\hat X^\an$ be an analytic open subset isomorphic to
    an open annulus.  Then $\hat X^\an\setminus A$ is an affinoid
    domain in $\hat X^\an$.
  \end{enumerate}
\end{lem}

\pf First we establish~(1).  Let us fix an isomorphism $B\cong\B(1)_+$. 
By~\cite[Lemma~3.5(c)]{bosch_lutk:uniformization}, for any 
$a\in K^\times$ with $|a| < 1$ the compact set
$\hat X^\an\setminus\bB(a)_+$ is an affinoid domain in $\hat X^\an$.  The
limit boundary of $\hat X^\an\setminus\bB(a)_+$ in $\hat X^\an$ is the Gauss
point $\|\cdot\|_{\val(a)}$ of $\bB(a)$; this coincides with the Shilov
boundary of $\hat X^\an\setminus\bB(a)_+$
by~\cite[Proposition~2.1.12]{thuillier:thesis}. 
The proof of Lemma~\ref{lem:skeleton.is.graph} shows that 
$\del_{\lim}(B) = \{x\}$ where $x = \lim_{r\to 0} \|\cdot\|_r$.

By the Riemann-Roch theorem, there exists a
meromorphic function on $\hat X$ which is regular away from $0\in\B(1)_+$ 
and which has a zero outside of $\B(1)_+$.  Fix such a function $f$, and
scale it so that $|f(x)| = 1$. 
By Corollary~\ref{cor:decreasing.logf}, the function $F(y) = -\log|f(y)|$
is a monotonically decreasing
function on $\Sigma(\bS(0)_+)\cong(0,\infty)$ such that
$\lim_{r\to 0} F(\|\cdot\|_r) = 0$.
The meromorphic function $f$ defines a finite
morphism $\phi: \hat X\to\bP^1$, which analytifies to a finite morphism
$\phi^\an: \hat X^\an\to \bP^{1,\an}$.  Let
$Y = \{ y\in \hat X^\an~:~|f(y)|\leq 1 \}$
be the inverse image of $\B(1)\subset \bP^{1,\an}$ under
$\phi^\an$, so $Y$ is an affinoid domain in $\hat X^\an$.
For $a\in\fm_R\setminus\{0\}$ the point $\|\cdot\|_{\val(a)}$ is the
Shilov boundary of $\hat X^\an\setminus\B(a)_+$, so 
$|f|\leq\|f\|_{\val(a)}$ on $\hat X^\an\setminus\bB(a)_+$.  Since
$\hat X^\an\setminus B\subset\hat X^\an\setminus\B(a)_+$ for all 
$a\in\fm_R\setminus\{0\}$ we have $|f|\leq \lim_{r\to 0}\|f\|_r = 1$ on 
$\hat X^\an\setminus B$.  Therefore 
$\hat X^\an\setminus B\subset Y$.

We claim that $\hat X^\an\setminus B$ is a connected component of $Y$.
Clearly it is closed in $Y$.  Since $f$ has finitely many zeros in $B$,
there exists $a\in\fm_R\setminus\{0\}$ such that $f$ is a unit on 
$\bS(a)_+\subset\B(1)_+$.  By Lemma~\ref{lem:unit.logf.constant} 
we have that $|f|> 1$ on $\bS(a)_+$, so 
$\hat X^\an\setminus B = (\hat X^\an\setminus\B(a))\cap Y$ is open in
$Y$.  Hence $\hat X^\an\setminus\B(1)_+$ is affinoid, being a connected
component of the affinoid domain $Y$.

We will reduce the second assertion to the first by doing surgery on
$\hat X^\an$, following the proof
of~\cite[Proposition~3.6.1]{berkovich:etalecohomology}.  
Let $A_1$ be a closed annulus inside of $A$, so 
$A\setminus A_1\cong\bS(a)_+\djunion\bS(b)_+$ for 
$a,b\in\fm_R\setminus\{0\}$.   Let $(X')^\an$ be the
analytic curve obtained by gluing $\hat X^\an\setminus A_1$ to two copies of
$\B(1)_+$ along the inclusions $\bS(a)_+\inject\B(1)_+$ and
$\bS(b)_+\inject\B(1)_+$.  One verifies easily that $(X')^\an$ is proper in
the sense of~\cite[\S3]{berkovich:analytic_geometry}, so $(X')^\an$ is the
analytification of a unique algebraic curve $X'$.  By construction
$\hat X^\an\setminus A$ is identified with the affinoid domain
$(X')^\an\setminus(\B(1)_+\djunion\B(1)_+)$ in $(X')^\an$, so we can apply~(1)
twice to $(X')^\an$ to obtain the result.\qed

\begin{rem} \label{rem:boundary.is.type2}
  Let $\sU$ be an affinoid domain in $\hat X^\an$ and let $x$ be a Shilov
  boundary point of $\sU$.  Since $\td\sH(x)$, the residue field of the
  completed residue field $\sH(x)$ at $x$, is isomorphic to the
  function field of an irreducible component of the canonical reduction of
  $\sU$, the point $x$ has type $2$.  Hence Lemma~\ref{lem:delete.open.affinoid}
  implies that if $A\subset\hat X^\an$ is an open ball or an open annulus
  then $\del_{\lim}(A)$ consists of either one or two type-$2$ points of
  $\hat X^\an$ since $\del_{\lim}(A)$ is the Shilov boundary of 
  $\hat X^\an\setminus A$.
\end{rem}

\smallskip
Recall that if $V$ is a semistable vertex set of $\hat X$ then there is a
retraction $\tau_V = \tau_{\Sigma(\hat X,V)} :\hat X^\an\to\Sigma(\hat X,V)$.

\begin{lem} \label{lem:exists.balls}
  Let $V$ be a semistable vertex set of $\hat X$ and let
  $x\in V$.  Then there are infinitely many open balls in the semistable
  decomposition for $\hat X$ which retract to $x$.
\end{lem}

\pf Suppose that there is at least one edge of $\Sigma(\hat X,V)$.  Deleting
all of the open annuli in the semistable decomposition of $\hat X$ 
yields an affinoid domain $Y$ by Lemma~\ref{lem:delete.open.affinoid}.  The
set $\tau_V\inv(x)$ is a connected component of $Y$, so $\tau_V\inv(x)$ is an
affinoid domain as well.  The Shilov boundary of $\tau_V\inv(x)$ agrees with
its limit boundary $\{x\}$ in $\hat X^\an$; by construction
$\tau_V\inv(x)\setminus\{x\}$ is a disjoint union of open balls, which are
the formal fibers of the canonical model of $\tau_V\inv(x)$
by~\cite[Lemme~2.1.13]{thuillier:thesis}.  Any nonempty curve over $k$ has
infinitely many points, so $\tau_V\inv(x)\setminus\{x\}$ is a disjoint union
of infinitely many open balls.

If $\Sigma(\hat X,V)$ has no edges then $\hat X^\an\setminus\{x\}$ is a
disjoint union of open balls.  Deleting one of these balls yields an
affinoid domain by Lemma~\ref{lem:delete.open.affinoid}, and the above
argument goes through.\qed

\paragraph \pf[of Theorem~\ref{thm:construct.models}]
First we prove that $\fX\mapsto V(\fX)$ is surjective, i.e.,\ that any
semistable vertex set comes from a semistable formal model.
Let $V$ be a semistable vertex set of $\hat X$, let 
$\Sigma = \Sigma(\hat X,V)$, and let 
$\tau = \tau_\Sigma:\hat X^\an\to\Sigma$ be the retraction.

\subparagraph[Case 1]   \label{par:construct.models.case1}
Suppose that $\Sigma$ has at least two edges.  Let $e$ be an edge in
$\Sigma$, let $A_0,A_1,\ldots,A_r$ ($r\geq 1$) be the open annuli in the
semistable decomposition of $\hat X$, and suppose that $\Sigma(A_0)$ is the interior
of $e$.  Then $\hat X\setminus(\bigcup_{i=1}^r A_i)$ is an affinoid
domain by Lemma~\ref{lem:delete.open.affinoid}, and $\tau\inv(e)$ is a
connected component of $\hat X\setminus(\bigcup_{i=1}^r A_i)$.
Hence $\tau\inv(e)$ is an affinoid domain in $\hat X^\an$.  Let
$\fY$ be its canonical model.  Let $x,y\in\hat X^\an$ be the endpoints of
$e$, so $\{x,y\} = \del_{\lim}(\tau\inv(e))$ is the Shilov boundary of
$\tau\inv(e)$, and $\tau\inv(e)\setminus\{x,y\}$ is a disjoint union of
open balls and the open annulus $A_0$.  By
\cite[Lemme~2.1.13]{thuillier:thesis}, the formal fibers of 
$\tau\inv(e)\to\bar\fY$ are the connected components of 
$\tau\inv(e)\setminus\{x,y\}$, so $\bar\fY$ has either one or two
irreducible components (depending on whether $x=y$) which intersect along
a single ordinary double point $\xi$
by Theorem~\ref{thm:inverse.image.reduction}.  Let
$\fC_x$ (resp.\ $\fC_y$) be the irreducible component of $\fY$
whose generic point is the reduction of $x$ (resp.\ $y$).  
Using the anti-continuity of the reduction map one sees that 
$\red\inv(\bar\fC_x\setminus\{\xi\}) = \tau\inv(x)$ and
$\red\inv(\bar\fC_y\setminus\{\xi\}) = \tau\inv(y)$.  It follows that the
formal affine subset $\fC_x\setminus\{\xi\}$ (resp.\ 
$\fC_y\setminus\{\xi\}$)  is the canonical model of the affinoid domain
$\tau\inv(x)$ (resp.\ $\tau\inv(y)$).  

Applying the above for every edge $e$ of $\Sigma$ allows us to glue the
canonical models of the affinoid domains $\tau\inv(e)$ together along the
canonical models of the affinoid domains $\tau\inv(x)$ corresponding to
the vertices $x$ of $\Sigma$.  Thus we obtain a semistable formal model
$\fX$ of $\hat X$ such that $V(\fX) = V$
(cf.\ Remark~\ref{rem:remarks.semistable}(1)). 

\subparagraph[Case 2]
Suppose that $\Sigma$ has one edge $e$ and two vertices
$x,y$.  Let $B_x,B_x'$ (resp.\ $B_y,B_y'$) be distinct open balls in the
semistable decomposition of $\hat X$ retracting to $x$ (resp.\ $y$), so 
$Y\coloneq\hat X^\an\setminus(B_x\cup B_y)$ and
$Y'\coloneq\hat X^\an\setminus(B_x'\cup B_y')$
are affinoid domains by Lemma~\ref{lem:delete.open.affinoid}.  
Let $\fY$ (resp.\ $\fY'$) be the canonical model of $Y$ (resp.\ $Y'$). 
Arguing as in Case~1 above, $\bar\fY$ and $\bar\fY{}'$ are affine curves
with two irreducible components intersecting along a single ordinary
double point $\xi$.  Furthermore, $Z = Y\cap Y'$ is an affinoid domain
whose canonical model $\fZ$ is obtained from $\fY$ (resp.\ $\fY'$) by
deleting one smooth point from each component.  Gluing $\fY$ to $\fY'$
along $\fZ$ yields the desired semistable formal model $\fX$ of $\hat X$.

\subparagraph[Case 3]
Suppose that $\Sigma$ has just one vertex $x$.
Let $B,B'$ be distinct open balls in the semistable decomposition of
$\hat X$, let $Y = \hat X^\an\setminus B$, let
$Y' = \hat X^\an\setminus B'$, and let $Z = Y\cap Y'$.  Gluing the
canonical models of $Y$ and $Y'$ along the canonical model of $Z$ gives us
our semistable formal model as in Case~2.

\subparagraph 
A semistable formal model of $\hat X$ is determined by its formal 
fibers~\cite[Lemma~3.10]{bosch_lutk:uniformization}, so
$\fX\mapsto V(\fX)$ is bijective.  It remains to prove that
$\fX$ dominates $\fX'$ if and only if $V(\fX')\subset V(\fX)$.  
If $\fX$ dominates $\fX'$ then $V(\fX')\subset V(\fX)$ by the
surjectivity and functoriality of the reduction map.  Conversely
let $V,V'$ be semistable vertex sets of $\hat X$ such that 
$V'\subset V$.  The corresponding semistable formal models 
$\fX,\fX'$ were constructed above by finding coverings $\sU$, $\sU'$ of 
$\hat X^\an$ by affinoid domains whose canonical models glue along the
canonical models of their intersections.  (Such a covering is called a 
\emph{formal covering} in~\cite{bosch_lutk:uniformization}.)  It is clear
that if $\sU$ refines $\sU'$, in the sense that every affinoid in $\sU$ is
contained in an affinoid in $\sU'$, then we obtain a morphism
$\fX\to\fX'$ of semistable formal models.  Therefore it suffices to show
that we can choose $\sU,\sU'$ such that $\sU$ refines $\sU'$ when
$V'\subset V$ in all of the cases treated above.  We will carry out this
procedure in the situation of Case~1, when $V$ is the union of $V'$ with a
type-$2$ point $x\in\Sigma'=\Sigma(\hat X,V')$ not contained in $V'$; the
other cases are similar and are left to the reader (cf.\ the proof of
Proposition~\ref{prop:relations.vertex.sets}). 

In the situation of Case~1, the formal covering corresponding to $V'$ is
the set 
\[ \sU' = \{\tau\inv(e)~:~e\text{ is an edge of }\Sigma'\}. \]
By Proposition~\ref{prop:relations.vertex.sets}(2) the skeleton
$\Sigma = \Sigma(\hat X, V)$ is a refinement of $\Sigma'$, obtained by
subdividing the edge $e_0$ containing $x$ to allow $x$ as a vertex.  Let
$e_1,e_2$ be the edges of $\Sigma$ containing $x$.  Then
$\tau\inv(e_1),\tau\inv(e_2)$ are affinoid domains in $\hat X^\an$
contained in $\tau\inv(e_0)$, so the formal covering
$\sU = \{ \tau\inv(e)~:~e\text{ is an edge of }\Sigma\}$ is a refinement
of $\sU'$, as desired.\qed

\paragraph[Stable models and the minimal skeleton]
Here we explain when and in what sense there exists a minimal semistable
vertex set of $X$.  Of course this question essentially reduces to the
existence of a stable model of $X$ when $X = \hat X$;
using~\cite{bosch_lutk:uniformization} we can also treat the case when $X$
is not proper.

\begin{defn*}
  Let $x\in X^\an$ be a type-$2$ point.  The \emph{genus} of $x$, denoted
  $g(x)$, is defined to be the genus of the smooth proper connected $k$-curve with
  function field $\td\sH(x)$, the residue field of the completed residue
  field $\sH(x)$ at $x$.
\end{defn*}

\begin{rem}
  Let $V$ be a semistable vertex set of $\hat X$ and let
  $x\in\hat X^\an$ be a type-$2$ point with positive genus.  Then
  $x\in V$, since otherwise $x$ admits a neighborhood which is isomorphic to
  an analytic domain in $\bP^{1,\an}$
  and the genus of any type-$2$ point in $\bP^{1,\an}$ is zero.
\end{rem}

\begin{rem} \label{rem:genus.formula}
  Let $\fX$ be a semistable formal model for $\hat X$, let
  $x\in V(\fX)$, and let $\bar\fC\subset\bar\fX$ be the irreducible
  component with generic point $\zeta = \red(x)$.  Then
  $\td\sH(x)$ is isomorphic to $\sO_{\bar\fX,\zeta}$
  by~\cite[Proposition~2.4.4]{berkovich:analytic_geometry}, so $g(x)$ is the
  genus of the normalization of $\bar\fC$.  It follows
  from~\cite[Theorem~4.6]{bosch_lutk:uniformization} that
  \begin{equation}\label{eq:genus.formula}
    g(\hat X) = \sum_{x\in V(\fX)} g(x) + g(\Sigma(\hat X,V))
  \end{equation}
  where $g(\hat X)$ is the genus of $\hat X$ and 
  $g(\Sigma(\hat X, V)) = \rank_\Z(H_1(\Sigma(\hat X,V),\Z))$ is the genus
  of $\Sigma(\hat X, V)$ as a 
  topological space (otherwise known as the cyclomatic number of the graph
  $\Sigma(\hat X,V)$).  The important
  equation~\eqref{eq:genus.formula} is known as the \emph{genus formula}. 
\end{rem}

\begin{defn}
  The \emph{Euler characteristic of $X$} is defined to be
  \[ \chi(X) = 2 - 2g(\hat X) - \#D. \]
\end{defn}

\begin{defn} \label{defn:stable.vertex.set}
  A semistable vertex set $V$ of $X$ is \emph{stable} if there
  is no $x\in V$ of genus zero and valence less than three in $\Sigma(X,V)$.
  We call the corresponding semistable decomposition of $X$ \emph{stable}
  as well.  A semistable formal model $\fX$ of $\hat X$ such that 
  $V(\fX)$ is a stable vertex set of $\hat X$ is called a 
  \emph{stable formal model}.
\end{defn}

A semistable vertex set $V$ of $X$ is 
\emph{minimal} if $V$ does not properly contain a semistable vertex set
$V'$. Any semistable vertex set contains a minimal one.

\begin{prop} \label{prop:delete.vertices}
  Let $V$ be a semistable vertex set of $X$ and let $x\in V$ be a point of
  genus zero.
  \begin{enumerate}
  \item Suppose that $x$ has valence one in $\hat\Sigma(X,V)$, let  
    $e$ be the edge adjoining $x$, and let $y$ be the other endpoint of
    $e$.  If $y\notin D$ 
    then $V\setminus\{x\}$ is a semistable vertex set of $X$ and 
    $\hat\Sigma(X,V\setminus\{x\})$ is the graph obtained from
    $\hat\Sigma(X,V)$ by removing $x$ and the interior of $e$.
  \item Suppose that $x$ has valence two in $\hat\Sigma(X,V)$, let
    $e_1,e_2$ be the edges adjoining $x$, and let $x_1$ (resp.\ $x_2$) be
    the other endpoint of $e_1$ (resp.\ $e_2$).  If
    $\{x_1,x_2\}\not\subset D$ then $V\setminus\{x\}$ is a semistable
    vertex set of $X$ and $\hat\Sigma(X,V\setminus\{x\})$ is the graph
    obtained from $\hat\Sigma(X,V)$ by joining $e_1,e_2$ into a single edge.
  \end{enumerate}
\end{prop}

\pf This is essentially~\cite[Lemma~6.1]{bosch_lutk:uniformization}
translated into our language.\qed

\smallskip
By a \emph{topological vertex} of a finite connected graph $\Gamma$ we
mean a vertex of valence at least $3$.   The set of topological vertices
only depends on the topological realization of $\Gamma$.

\begin{thm}[Stable reduction theorem] \label{thm:stable.reduction}
  There exists a semistable vertex set of $X$.  If $V$ is a minimal
  semistable vertex set of $X$ then:
  \begin{enumerate}
  \item If $\chi(X)\leq 0$ then $\Sigma(X,V)$ is the set of points in
    $X^\an$ that do not admit an affinoid neighborhood isomorphic to
    $\B(1)$.
  \item If $\chi(X) < 0$ then $V$ is stable and 
    \[ V = \{ x\in\Sigma(X,V)~:~x\text{ is a topological vertex of }
    \Sigma(X,V)\text{ or } g(x) > 0 \}. \]
  \end{enumerate}
\end{thm}

\begin{cor} \label{cor:minimal.skeleton}
  If $\chi(X)\leq 0$ then there is a unique set-theoretic minimal skeleton
  of $X$, and if $\chi(X) < 0$ then there is a unique stable vertex set of
  $X$. 
\end{cor}

\smallskip
\pf[of Theorem~\ref{thm:stable.reduction}]
The existence of a semistable vertex set of $\hat X$ follows from the
classical theorem of Deligne and
Mumford~\cite{deligne_mumford:irreducibility} as  
proved analytically (over a non-noetherian rank-$1$ valuation ring)
in~\cite[Theorem~7.1]{bosch_lutk:uniformization}.  
The existence of a semistable vertex set of $X$ then follows
from Proposition~\ref{prop:relations.vertex.sets}(3).  
Let $V$ be a minimal semistable vertex set of $X$ and let 
$\Sigma = \Sigma(X,V)$.  If $\chi(X) < 0$ then one
applies Proposition~\ref{prop:delete.vertices} in the standard way to prove the
second assertion, and if $\chi(X)\leq 0$
then Proposition~\ref{prop:delete.vertices}(1) guarantees that every genus-zero
vertex of $\Sigma$ has valence at least two.

Suppose that $\chi(X)\leq 0$.  
Let $\Sigma'$ be the set of points of
$X^\an$ that do not admit an affinoid neighborhood isomorphic to $\B(1)$.
By Lemma~\ref{lem:skeleton.nice}(4) we have $\Sigma'\subset\Sigma$.  Let
$x\in\Sigma$, and suppose that $x$ admits an affinoid neighborhood $U$
isomorphic to $\B(1)$.  We will show by way of contradiction that $\Sigma$
has a vertex of valence less than two in $U$ (any vertex contained in $U$
has genus zero); in fact we will show that 
$\Sigma\cap U$ is a tree.  Let $y$ be the Gauss point of $U$.  If 
$y\in\Sigma$ then we may replace $V$ by $V\cup\{y\}$
by Proposition~\ref{prop:relations.vertex.sets}(2) to assume that $y\in V$. 
Since $U$ is closed and any connected component of $X^\an\setminus V$ that
intersects $U$ is contained in $U$, the retraction $\tau_\Sigma:X^\an\to\Sigma$
restricts to a retraction $U\to U\cap\Sigma$.  Since $U$ is contractible, 
$U\cap\Sigma$ is a tree as claimed.\qed

\begin{rem} \label{rem:tate.curves}
  If $\chi(X) = 0$ then either $g(X) = 0$ and $\#D = 2$ or
  $g(X) = 1$ and $\#D = 0$.  In the first case, the skeleton of 
  $X\cong\G_m$ is the line connecting $0$ and $\infty$, and any type-$2$
  point on this line is a minimal semistable vertex set.  In the second case,
  $X = \hat X$ is an elliptic curve with respect to some choice of
  distinguished point $0\in X(K)$.  If $X$ has good reduction then
  there is a unique point $x\in X^{\an}$ with $g(x) = 1$; in this case $\{x\}$
  is the unique stable vertex set of $X$ and $\Sigma(X,\{x\}) = \{x\}$.

  Suppose now that $(X,0)$ is an elliptic curve with multiplicative
  reduction, i.e.,\ $X$ is a Tate curve.  By Tate's uniformization
  theory~\cite[\S9.7]{bgr:nonarch}, 
  there is a unique $q = q_X\in K^\times$ with $\val(q) > 0$ and an \'etale
  morphism $u: \G_m^\an\to X^\an$ which is a 
  homomorphism of group objects (in the category of $K$-analytic spaces)
  with kernel $u\inv(0) = q^\Z$.  For brevity we will often write
  $X^\an\cong\G_m^\an/q^\Z$.  The so-called Tate parameter $q$ 
  is related to the $j$-invariant $j=j_X$ of $X$ in such a way that
  $\val(q) = -\val(j)$ (it is the $q$-expansion of the modular function
  $j$).  Let   $Z$ be the retraction of the set $q^\Z$ onto 
  the skeleton of $\G_m$, i.e.,\ the collection of Gauss points of the
  balls $\B(q^n)$ for $n\in\Z$.  Then $\G_m^\an\setminus Z$ is
  the disjoint union of the open annuli $\{\bS(q^{n+1},q^n)_+\}_{n\in\Z}$
  and infinitely many open balls, and every connected component of
  $\G_m^\an\setminus Z$ maps
  isomorphically onto its image in $X^\an$.  It follows that 
  $X^\an\setminus\{u(1)\}$ is a disjoint union of an open annulus $A$
  isomorphic to $\bS(q)_+$ and infinitely many open balls.  Hence
  $V = \{u(1)\}$ is a (minimal) semistable vertex set of $X$, and the
  associated (minimal) skeleton $\Sigma$ is a circle of circumference 
  $\val(q) = -\val(j_E)$.  We have $u(1) = \tau_\Sigma(0)$, so any
  type-$2$ point on $\Sigma$ is a minimal semistable vertex set, as any
  such point is the retraction of a $K$-point of $X$ (which we could have
  chosen to be $0$).

  See also~\cite[Example~7.20]{gubler:lchs}.
\end{rem}

\begin{rem} \label{rem:Mgtrop} Given a smooth complete curve $\hat{X}/K$
  of genus $g$ and a subset $D$ of `marked points' of $\hat X(K)$ 
  satisfying the inequality $2-2g-n \leq 0$, where $n = \#D \geq 0$, one
  obtains a canonical pair $(\Gamma, w)$ consisting of an abstract metric
  graph and a vertex weight function, where $\Gamma = \Sigma(\hat{X}
  \setminus D,V)$ is the minimal skeleton of $\hat X\setminus D$ and $w :
  \Gamma \to \ZZ_{\ge 0}$ takes $x \in \Gamma$ to $0$ if $x \not\in V$ and
  to $g(x)$ if $x \in V$.  (A closely related construction can be found in
  \cite[\S{2}]{tyomkin:correspondence_thms}.)  If $2-2g-n < 0$, this gives a canonical `abstract
  tropicalization map' $\trop : M_{g,n} \to M_{g,n}^{\rm trop}$, where $M_{g,n}^{\rm
    trop}$ is the moduli space of $n$-pointed tropical curves of genus $g$
  as defined, for example, in
  \cite[\S{3}]{caporaso:tropical_moduli_spaces}.  The map $\trop : M_{g,n} \to
  M_{g,n}^{\rm trop}$ is certainly deserving of further study.
\end{rem}

\paragraph[Application to the local structure theory of $X$]
The semistable reduction theorem and its translation into the language of
semistable vertex sets yields the following information about the local
structure theory of an analytic curve.
(Conversely, one can study the local structure of an analytic curve
directly and derive the semistable reduction theorem:
see~\cite{temkin:stable_modification}.)

\begin{cor} \label{cor:fundamental.nbhds}
  Let $x\in X^\an$.  There is a fundamental system of open neighborhoods
  $\{U_\alpha\}$ of $x$ of the following form:
  \begin{enumerate}
  \item If $x$ is a type-$1$ or a type-$4$ point then the $U_\alpha$ are
    open balls.
  \item If $x$ is a type-$3$ point then the $U_\alpha$ are open annuli
    with $x\in\Sigma(U_\alpha)$.
  \item If $x$ is a type-$2$ point then $U_\alpha = \tau_V\inv(W_\alpha)$
    where $W_\alpha$ is a simply-connected open neighborhood of $x$ in
    $\Sigma(X,V)$ for some semistable vertex set $V$ of $X$ containing
    $x$, and each $U_\alpha\setminus\{x\}$ is a
    disjoint union of open balls and open annuli. 
  \end{enumerate}
\end{cor}

\pf Since $X$ has a semistable decomposition, if $x$ is a point of type
$1$, $3$, or $4$ then $x$ has a neighborhood isomorphic to an open annulus
or an open ball.  Hence we may assume that $X = \bP^1$ and $x\in\B(1)_+$.
By~\cite[Proposition~1.6]{baker_rumely:book} the set of open balls with
finitely many closed balls removed forms a basis for the topology on
$\B(1)_+$; assertions~(1) and~(2) follow easily from this.  

Let $f$ be a meromorphic function on $X$; deleting the zeros and poles of
$f$, we may assume that $f$ is a unit on 
$X$.  Let $F = \log |f|: X^\an\to\R$ and let $U = F\inv((a,b))$ for some
interval $(a,b)\subset\R$.  Let $x$ be a type-$2$ point contained in $U$.
Since such $U$ form a sub-basis for the topology on $X^\an$ it suffices to
prove that there is a neighborhood of $x$ of the form described in~(3)
contained in $U$.  Let $V$ be a semistable vertex set for $X$ containing
$x$.  By Proposition~\ref{prop:annulus.skeleton}
and Lemma~\ref{lem:unit.logf.constant} we have that $F$ is affine-linear on the
edges of $\Sigma(X,V)$ and that $F$ factors through
$\tau_V:X^\an\to\Sigma(X,V)$.  Therefore if $W$ is any simply-connected
neighborhood of $x$ in $\Sigma(X,V)$ contained in $U = F\inv((a,b))$ then
$\tau_V\inv(W)\subset U$.  If we assume in addition that the intersection of
$W$ with any edge of $\Sigma$ adjoining $x$ is a half-open interval with
endpoints in $G$ then $\tau_V\inv(W)\setminus\{x\}$ is a disjoint union of
open balls and open annuli.\qed

\begin{defn} \label{defn:simple.nbhd}
  A neighborhood of $x\in X^\an$ of the form described
  in Corollary~\ref{cor:fundamental.nbhds} is called a 
  \emph{simple neighborhood}  of $x$.
\end{defn}

\paragraph \label{par:nbhd.type2}
A simple neighborhood of a type-$2$ point $x\in X^\an$ has
the following alternative description.  Let $V$ be a semistable vertex set
containing $x$ and let $W$ be a simply-connected neighborhood of $x$ in
$\Sigma(X,V)$ such that the intersection of $W$ with any edge adjoining
$x$ is a half-open interval with endpoints in $G$, so $U = \tau_V\inv(W)$ is
a simple neighborhood of $x$.  Adding the boundary of $W$ to $V$, we may
assume that the connected components of $U\setminus\{x\}$ are connected
components of $X^\an\setminus V$.  Let $\fX$ be the semistable formal
model of $\hat X$ associated to $V$ and let $\bar\fC\subset\bar\fX$ be the
irreducible component with generic point $\red(x)$.  Since $W$ contains no
loop edges of $\Sigma(X,V)$, the component $\bar\fC$ is smooth.  The connected
components of $\hat X^\an\setminus V$ are the formal fibers of $\fX$, so
it follows from the anti-continuity of $\red$ that 
$U = \red\inv(\bar\fC)$ and that $\pi_0(U\setminus\{x\})\isom\bar\fC(k)$.
To summarize:

\begin{lem*}
  A simple neighborhood $U$ of a type-$2$ point $x\in X^\an$ is the
  inverse image of a smooth irreducible component $\bar\fC$ of the special
  fiber of a semistable formal model $\fX$ of $\hat X$.  Furthermore, we
  have   $\pi_0(U\setminus\{x\})\isom\bar\fC(k)$.
\end{lem*}

\section{The metric structure on an analytic curve}
\label{sec:metric.structure}

The set of all skeleta $\{\Sigma(X,V)\}_V$ is a filtered directed system
under inclusion by Proposition~\ref{prop:relations.vertex.sets}. 
For $U$ a one-dimensional $K$-analytic space, define the set of
\emph{skeletal points} $\HH_\circ(U)$ of $U$ to be
the set of points of $U$ of types $2$ and $3$, and the set of
\emph{norm points} to be $\HH(U)\coloneq U\setminus U(K)$.  
When $U=X$
the latter are the points that arise from norms on the function field 
$K(X)$ which extend the given absolute value on $K$, and the
following corollary explains the former terminology:

\begin{cor} \label{cor:direct.limit.skeleta}
  We have
  \[ \HH_\circ(X^\an) = \bigcup_V \Sigma(X,V) = \varinjlim_V\Sigma(X,V) \]
  as sets, where $V$ runs over all semistable vertex sets of $X$.
\end{cor}

\pf Any point of $\Sigma(X,V)$ has type $2$ or $3$, and any type-$2$ point
is contained in a semistable vertex
set by Proposition~\ref{prop:relations.vertex.sets}(3).  Let $x$ be a type-$3$
point.  Then $x$ is contained in an open ball or an open annulus in a
semistable decomposition of $X^\an$.  The semistable decomposition can
then be refined as in the proof
of Proposition~\ref{prop:relations.vertex.sets}(3) to produce a skeleton that
includes $x$.\qed

\smallskip
By Proposition~\ref{prop:relations.vertex.sets}(1), the set of all skeleta
$\{\Sigma(X,V)\}_V$ is also an inverse system with respect to the natural
retraction maps.  Although not logically necessary for anything else in
this paper, the following folklore counterpart to
Corollary~\ref{cor:direct.limit.skeleta} is conceptually important.  For a
higher-dimensional analogue (without proof) in the case ${\rm char}(K)=0$, see
\cite[Appendix A]{kontsevich_soibelman:affinestructures}, 
and see~\cite[Corollary~3.2]{boucksom_favre_johnsson:metrics} in general.
See also~\cite{hrushovsky_loeser:tame_topology}.

\begin{thm} \label{thm:inverse.limit.skeleta}
The natural map
\[ u : \hat X^\an \to \varprojlim_V\Sigma(\hat X,V) \] is a homeomorphism
of topological spaces, where $V$ runs over all semistable vertex sets of
$\hat X$.
\end{thm}

\pf The map $u$ exists and is continuous by the universal property of
inverse limits.  It is injective because given any two points $x \neq y$
in $\hat X^\an$, one sees easily that there is a semistable vertex set $V$ such
that $x$ and $y$ retract to different points of $\Sigma(\hat X,V)$.  Since
$\hat X^\an$ is compact and each 
individual retraction map $\hat X^\an \to \Sigma(\hat X,V)$ is continuous and
surjective, it follows from 
\cite[\S{9.6} Corollary~2]{bourbaki:generaltopology} 
that $u$ is also surjective.  By Proposition 8 in \S{9.6} of loc.~cit., 
the space $\varprojlim_V\Sigma(\hat X,V)$ is compact.
Therefore $u$ is a continuous bijection between compact (Hausdorff) spaces, hence a homeomorphism
(cf. Corollary 2 in \S{9.4} of loc.~cit.).
\qed

\smallskip
\paragraphtoc[The metric structure on $\HH_\circ(X^\an)$]
Let $V\subset V'$ be semistable vertex sets of $X$.
By Proposition~\ref{prop:relations.vertex.sets}(3) every edge $e$ of
$\Sigma(X,V)$ includes isometrically into an edge of $\Sigma(X,V')$.  Let
$x,y\in\Sigma(X,V)$ and let $[x,y]$ be a shortest path
from $x$ to $y$ in $\Sigma(X,V)$.  Then $[x,y]$ is also a shortest path
in $\Sigma(X,V')$: if there were a shorter path $[x,y]'$ in $\Sigma(X,V')$
then $[x,y]\cup[x,y]'$ would represent a homology class in 
$H_1(\Sigma(X,V'),\Z)$ that did not exist in $H_1(\Sigma(X,V),\Z)$, which
is impossible by the genus formula~\eqref{eq:genus.formula}.  Therefore
the inclusion $\Sigma(X,V)\inject\Sigma(X,V')$ is an isometry (with
respect to the shortest-path metrics), so
by Corollary~\ref{cor:direct.limit.skeleta} we obtain a natural metric 
$\rho$ on $\HH_{\circ}(X^\an)$, called the \emph{skeletal metric}.

Let $V$ be a semistable vertex set and let $\tau = \tau_V: X^\an\to\Sigma(X,V)$ be
the retraction onto the skeleton.  If  $x,y\in \HH_\circ(X^\an)$ are not
contained in the same connected component of $X^\an\setminus\Sigma(X,V)$
then a shortest path from $x$ to $y$ in a larger skeleton must go through
$\Sigma(X,V)$.  It follows that
\begin{equation} \label{eq:metric.on.skeleton}
  \rho(x,y) = \rho(x,\tau(x)) + \rho(\tau(x),\tau(y)) +
  \rho(\tau(y),y). 
\end{equation}

\begin{rem}
  \begin{enumerate}
  \item By definition any skeleton includes isometrically into
    $\HH_\circ(X^\an)$. 
  \item It is important to note that the metric topology on
    $\HH_\circ(X^\an)$ is \emph{stronger} than the subspace topology. 
  \end{enumerate}
\end{rem}

We can describe the skeletal metric locally as follows.  By
Berkovich's classification theorem, any point $x\in \HH(\A^{1,\an})$
is a limit of Gauss points of balls of radii $r_i$ converging to
$r\in(0,\infty)$.  We define $\diam(x) = r$.  Any two points
$x \neq y\in \A^{1,\an}$ are contained in a unique smallest closed ball; its
Gauss point is denoted $x\vee y$.  For $x,y\in \HH(\A^{1,\an})$ we
define 
\[ \rho_p(x,y) = 2\log(\diam(x\vee y)) - \log(\diam(x)) - \log(\diam(y)). \]
Then $\rho_p$ is a metric on $\HH(\A^{1,\an})$, called the 
\emph{path distance metric}; see~\cite[\S2.7]{baker_rumely:book}.  If $A$ 
is a standard open ball or standard generalized open annulus then the restriction of
$\rho_p$ to $\HH(A)$ is called the \emph{path distance metric} on
$\HH(A)$. 

\begin{prop} \label{prop:restrict.metric}
  Let $A\subset X^\an$ be an analytic domain isomorphic to a standard open
  ball or a standard generalized open annulus.  Then the skeletal metric on
  $\HH_\circ(X^\an)$ and the path distance metric on $\HH(A)$
  restrict to the same metric on $\HH_\circ(A)$. 
\end{prop}

\pf Let $V$ be a semistable vertex set
containing the limit boundary of $A$
(cf.~Remark~\ref{rem:boundary.is.type2}).
Then $V\setminus(V\cap A)$ is a semistable vertex set since the connected
components of $A\setminus(V\cap A)$ are connected components of 
$X^\an\setminus V$.  Hence we may and do assume that $A$ is a connected
component of $X^\an\setminus V$.  Suppose that $A$ is an open ball, and
fix an isomorphism $A\cong\B(a)_+$.  Let $x,y\in A$ be type-$2$ points.  
\begin{enumerate}
\item Suppose that $x\vee y\in\{x,y\}$; without loss of generality we may
  assume that $x = x\vee y$.  After recentering, we may assume in addition
  that $x$ is the Gauss point of $\B(b)$ and that $y$ is the Gauss point of
  $\B(c)$.  Then the standard open annulus $A' = \bB(b)_+\setminus\bB(c)$
  is a connected component of $A\setminus\{x,y\}$, which breaks up into a
  disjoint union of open balls and the open annuli $A'$ and
  $\B(a)_+\setminus\B(b)$.  Hence $V\cup\{x,y\}$ 
  is a semistable vertex set, and $\Sigma(A')$ is the interior of the edge
  $e$ of $\Sigma(X,V\cup\{x,y\})$ with endpoints $x,y$.
  Therefore $\rho(x,y)$ is the logarithmic modulus of $A'$, which agrees
  with $\rho_p(x,y) = \log(\diam(x)) - \log(\diam(y))$.
\item Suppose that $z = x\vee y\notin\{x,y\}$.  Then 
  $A\setminus\{x,y,z\}$ is a disjoint union of open balls and three open
  annuli, two of which connect $x,z$ and $y,z$.  As above we have
  $\rho(x,y) = \rho(x,z) + \rho(y,z)$, which is the same as 
  $\rho_p(x,y) = (\log(\diam(z))-\log(\diam(x))) + (\log(\diam(z))-\log(\diam(y)))$.
\end{enumerate}
Since the type-$2$ points of $A$ are dense~\cite[Lemma~1.8]{baker_rumely:book},
this proves the claim when $A$ is an open ball in a semistable
decomposition of $X$.  The proof when $A$ is a generalized open annulus in
a semistable decomposition of $X$ has more cases but is not essentially
any different, so it is left to the reader.\qed

\smallskip
Since Proposition~\ref{prop:restrict.metric} did not depend on the choice of
isomorphism of $A$ with a standard generalized open annulus, we obtain:

\begin{cor}
  Any isomorphism of standard open balls or standard generalized open
  annuli induces an isometry with respect to the path distance metric.
\end{cor}

\smallskip
In particular, if $A$ is an (abstract) open ball or generalized open
annulus then we can speak of the path distance metric on $\HH(A)$.

\begin{cor} \label{cor:extend.metric}
  The metric $\rho$ on $\HH_\circ(X^\an)$ extends in a unique way to a
  metric on $\HH(X^\an)$.
\end{cor}

\pf Let $x,y\in \HH(X^\an)$ and let $V$ be a semistable vertex set
of $X$.  If $x,y$ are contained in the same connected component
$B\cong\B(1)_+$ of $X^\an\setminus\Sigma(X,V)$ then we set
$\rho(x,y) = \rho_p(x,y)$.  Otherwise we set
\[ \rho(x,y) = \rho_p(x,\tau_V(x)) + \rho(\tau_V(x),\tau_V(y)) +
\rho_p(\tau_V(y),y) \]
where we have extended the path distance metric $\rho_p$ on a connected
component $B$ of $X^\an\setminus\Sigma(X)$ to its closure $B\cup\tau_V(B)$
by continuity (compare the proof of Lemma~\ref{lem:skeleton.is.graph}).
By~\eqref{eq:metric.on.skeleton} and Proposition~\ref{prop:restrict.metric} this
function extends $\rho$.  We leave it to the reader to verify that $\rho$
is a metric on $\HH(X^\an)$.\qed

\paragraph
A \emph{geodesic segment from $x$ to $y$} in a metric space $T$ is the
image of an isometric embedding $[a,b]\inject T$ with $a\mapsto x$ and
$b\mapsto y$.  We often identify a geodesic segment with its image in $T$.
Recall that an \emph{$\R$-tree} is a metric space $T$ with the following
properties: 
\begin{enumerate}
\item For all $x,y\in T$ there is a unique geodesic segment $[x,y]$ from
  $x$ to $y$.
\item For all $x,y,z\in T$, if $[x,y]\cap[y,z] = \{y\}$ then
  $[x,z] = [x,y]\cup[y,z]$.
\end{enumerate} 
See~\cite[Appendix~B]{baker_rumely:book}.
It is proved in~\S1.4 of loc.\ cit.\ that $\HH(\B(1))$ is
an $\R$-tree under the path distance metric.  It is clear that any
path-connected subspace of an $\R$-tree is an $\R$-tree, so
if $A$ is an open ball or a generalized open annulus then 
$\HH(A)$ is an $\R$-tree as well.

\begin{prop} \label{prop:locally.Rtree}
  Every point $x\in \HH(X^\an)$ admits a fundamental system of
  simple neighborhoods $\{U_\alpha\}$ in $X^\an$ such that
  $U_\alpha\cap \HH(X^\an)$ is an $\R$-tree under the restriction
  of $\rho$. 
\end{prop}

The definition of a simple neighborhood of a point $x\in X^\an$ is found
in~\parref{defn:simple.nbhd}. 

\pf If $x$ has type $3$ or $4$ then a simple neighborhood of $x$ is an
open ball or an open annulus, so the proposition follows 
from Corollary~\ref{cor:fundamental.nbhds} and
Proposition~\ref{prop:restrict.metric}. 
Let $x$ be a type-$2$ point and let $V$ be a semistable vertex set of $X$
containing $x$.  For small enough $\epsilon > 0$ the set
$W = \{ y\in\Sigma(X,V)~:~\rho(x,y) < \epsilon \}$ is
simply-connected; fix such an $\epsilon\in G$, and let $U = \tau_V\inv(W)$.
Then $U$ is a simple neighborhood of $x$.  We claim that $\HH(U)$
is an $\R$-tree.  Any connected component $A$ of $U\setminus\{x\}$ is an
open ball or an open annulus, so $\HH(A)$ is an $\R$-tree.
Moreover $\HH(A)\cup\{x\}$ is isometric to a path-connected
subspace of $\HH(\B(1))$ as in the proof
of Lemma~\ref{lem:skeleton.is.graph}; it 
follows that $\HH(A)\cup\{x\}$ is an $\R$-tree.  Therefore
$\HH(U)$ is a 
collection of $\R$-trees joined together at the single point $x$, and the
hypotheses on $W$ along with~\eqref{eq:metric.on.skeleton}
imply that if $y,z\in \HH(U)$ are contained in
different 
components of $U\setminus\{x\}$ then $\rho(y,z) = \rho(y,x)+\rho(x,z)$.
It is clear that such an object is again an $\R$-tree.\qed

\begin{cor} \label{cor:geodesic.in.skeleton}
  Let $x,y\in \HH_\circ(X^\an)$ and let $\Sigma = \Sigma(X,V)$ be a
  skeleton containing $x$ and $y$.  Then any geodesic segment from $x$ to
  $y$ is contained in $\Sigma$.
\end{cor}

\pf Any path from $x$ to $y$ in $\Sigma$ is by definition a
geodesic segment.  If $x,y$ are contained in
an open subset $U$ such that $\HH(U)$ is an $\R$-tree then the path
from $x$ to $y$ in $\Sigma\cap U$ is the
unique geodesic segment from $x$ to $y$ in $\HH(U)$.  The general
case follows by covering a geodesic segment from $x$ to $y$ by (finitely many) such
$U$.\qed 

\paragraphtoc[Tangent directions and the Slope Formula]
\label{par:tangent.directions}
Let $x\in \HH(X^\an)$.  A 
\emph{nontrivial geodesic segment starting at $x$} is a geodesic segment 
$\alpha:[0,a]\inject \HH(X^\an)$ with $a > 0$ such that 
$\alpha(0) = x$.  We say that two nontrivial geodesic segments
$\alpha,\alpha'$ starting at $x$ are \emph{equivalent at $x$} if $\alpha$
and $\alpha'$ agree on a neighborhood of $0$.
Following~\cite[{\S}B.6]{baker_rumely:book}, we define the set of
\emph{tangent directions at $x$} to be the set $T_x$ of nontrivial
geodesic segments starting at $x$ up to equivalence at $x$. It is clear
that $T_x$ only depends on a neighborhood of $x$ in $X^\an$.

\begin{lem} \label{lem:tangent.directions}
  Let $x\in \HH(X^\an)$ and let $U$ be a simple neighborhood of $x$
  in $X^\an$.   Then $[x,y]\mapsto y$ establishes a bijection 
  $T_x\isom \pi_0(U\setminus\{x\})$.  Moreover,
  \begin{enumerate}
  \item If $x$ has type $4$ then there is only one tangent direction at
    $x$.
  \item If $x$ has type $3$ then there are two tangent directions at $x$.
  \item If $x$ has type $2$ then $U = \red\inv(\bar\fC)$ for a smooth
    irreducible component $\bar\fC$ of the special fiber of a semistable
    formal model $\fX$ of $\hat X$ by~\parref{par:nbhd.type2}, and 
    $T_x \isom \pi_0(U\setminus\{x\}) \isom \bar\fC(k)$.
  \end{enumerate}
\end{lem}

\pf We will assume for simplicity that $\HH(U)$ is an $\R$-tree
(i.e.,\ that the induced metric on $\HH(U)$ agrees with the
shortest-path metric);
the general case reduces to this because $U$ is contractible.
The bijection $T_x\isom\pi_0(\HH(U)\setminus\{x\})$ is proved
in~\cite[{\S}B.6]{baker_rumely:book}.  A connected component $B$ of 
$U\setminus\{x\}$ is an $\R$-tree by Proposition~1.13 of 
loc.\ cit.\ and the type-$1$ points of $B$ are leaves, so 
$\pi_0(\HH(U)\setminus\{x\}) = \pi_0(U\setminus\{x\})$.
Parts~(1) and~(2) are proved in \S1.4 of loc.\ cit, and
part~(3) is~\parref{par:nbhd.type2}.\qed

\paragraph\label{par:reduce.rational.func}
With the notation in Lemma~\ref{lem:tangent.directions}(3),
we have a canonical identification of $\td\sH(x)$ with the function field
of $\bar\fC$ by~\cite[Proposition~2.4.4]{berkovich:analytic_geometry}.  
Hence we have an identification $\xi\mapsto\ord_\xi$
of $\bar\fC(k)$ with the set $\DV(\td\sH(x)/k)$ of nontrivial discrete
valuations $\td\sH(x)\surject\Z$ inducing the trivial valuation on $k$.
One can prove that the composite bijection
$T_x\isom\DV(\td\sH(x)/k)$ is independent of the choice of $U$.  The
discrete valuation corresponding to a tangent direction $v\in T_x$ will be
denoted $\ord_v:\td\sH(x)\to\Z$. 

Let $x\in X^\an$ be a type-$2$ point and let $f$ be an analytic function
in a neighborhood of $x$.  Let $c\in K^\times$ be a scalar such that 
$|f(x)| = c$.  We define $\td f_x\in\td\sH(x)$ to be the residue of
$c\inv f$, so $\td f_x$ is only defined up to
multiplication by a nonzero scalar in $k$.  However if 
$\ord:\td\sH(x)\to\Z$ is a nontrivial discrete valuation trivial on $k$
then $\ord(\td f_x)$ is intrinsic to $f$.

\begin{defn}
  A function $F:X^\an\to\R$ is \emph{piecewise affine} provided that for
  any geodesic segment $\alpha:[a,b]\inject \HH(X^\an)$ the
  pullback $F\circ\alpha:[a,b]\to\R$ is piecewise affine.  The
  \emph{outgoing slope} of a piecewise affine function $F$ at a point
  $x\in \HH(X^\an)$ along a tangent direction $v\in T_x$ is defined
  to be 
  \[ d_vF(x) = \lim_{\epsilon\to 0} (F\circ\alpha)'(\epsilon) \]
  where $\alpha:[0,a]\inject X^\an$ is a nontrivial geodesic segment
  starting at $x$ which represents $v$.  We say that a piecewise affine
  function $F$ is \emph{harmonic} at a point $x\in X^\an$ provided that
  the outgoing slope $d_vF(x)$ is nonzero for
  only finitely many $v\in T_x$, and $\sum_{v\in T_x} d_vF(x) = 0$.
  We say that $F$ is \emph{harmonic} if it is harmonic for all 
  $x\in\HH(X^\an)$. 
\end{defn}

\begin{thm}[Slope Formula]  \label{thm:PL}
  Let $f$ be an algebraic function on $X$ with no zeros or poles and let
  $F = -\log |f|: X^\an\to\R$.  Let $V$ be a semistable vertex set of $X$
  and let $\Sigma = \Sigma(X,V)$.  Then:
  \begin{enumerate}
  \item $F = F\circ\tau_\Sigma$ where $\tau_\Sigma:X^\an\to\Sigma$ is the retraction.
  \item $F$ is piecewise affine with integer slopes, and $F$ is affine-linear on
    each edge of $\Sigma$. 
  \item If $x$ is a type-$2$ point of $X^\an$ and $v\in T_x$ then
    $d_v F(x) = \ord_v(\td f_x)$.
  \item $F$ is harmonic.
  \item Let $x\in D$, let $e$ be the ray in $\Sigma$ whose closure in
    $\hat X$ contains $x$, let $y\in V$ be the other endpoint of $e$, and
    let $v\in T_y$ be the tangent direction represented by $e$.  Then $d_v
    F(y) = \ord_x(f)$.
  \end{enumerate}
\end{thm}

\pf The first claim follows from Lemma~\ref{lem:unit.logf.constant} and the
fact that a unit on an open ball has constant absolute value.  The
linearity of $F$ on edges of $\Sigma$ is Proposition~\ref{prop:simple.PL}.
Since $F = F\circ\tau_\Sigma$ we have that $F$ is constant in a neighborhood of
any point of type $4$, and any geodesic segment contained in $\HH_\circ(X^\an)$
is contained in a skeleton by Corollary~\ref{cor:geodesic.in.skeleton}, so $F$ is
piecewise affine.  The last claim is Proposition~\ref{prop:simple.PL}(2).
The harmonicity of $F$ is proved as follows: if $x\in X^\an$ has
type $4$ then $x$ has one tangent direction and $F$ is locally constant in
a neighborhood of $x$, so $\sum_{v\in T_x} d_v F(x) = 0$.  If $x$ has type
$3$ then $x$ is contained in the interior of an edge $e$ of a skeleton,
and the two tangent directions $v,w$ at $x$ are represented by the two
paths emanating from $x$ in $e$; since $F$ is affine on $e$ we have
$d_v F(x) = -d_w F(x)$.  The harmonicity of $F$ at type-$2$ points is an
immediate consequence of~(3) and the fact that the divisor of a
meromorphic function on a smooth complete curve has degree zero.

The heart of this theorem is~(3), which again is essentially a result of
Bosch and L\"utkebohmert. Let $x$ be a type-$2$ point of
$X^\an$, let $U$ be a simple neighborhood of $x$, and let $\fX$ be a
semistable formal model of $\hat X$ such that $x\in V(\fX)$ and 
$U = \red\inv(\bar\fC)$ where $\bar\fC$ is the smooth irreducible
component of $\bar\fX$ with generic point $\red(x)$.  We may and do assume
that $V(\fX)$ is a semistable vertex set of $X$ containing $V$.  Let 
$\bar\fC{}'\subset\bar\fC$ be the affine curve obtained by deleting all
points $\xi\in\bar\fC$ which are not smooth in $\bar\fX$ and let 
$\fC'$ be the induced formal affine subscheme of $\fX$.  Then 
$(\fC')^\an = \red\inv(\bar\fC{}') = \tau_{V(\fX)}\inv(x)$ is an affinoid domain in
$X^\an$ with Shilov boundary $\{x\}$.  
If we scale $f$ such that $|f(x)| = 1$ then
  $f$ and $f\inv$ both have
supremum norm $1$ on $\tau_{V(\fX)}\inv(x)$.  It follows that the residue 
$\td f_x$ of $f$ is a unit on $\bar\fC{}'$, so 
$\ord_\zeta(\td f_x) = 0$ for all $\zeta\in\bar\fC{}'(k)$.   
By~(1) we have that $F$ is constant
on $\tau_{V(\fX)}\inv(x)$, so $d_v F(x) = \ord_v(\td f_x) = 0$ for all $v\in T_x$ 
corresponding to closed points of $\bar\fC{}'$.

Now let $v\in T_x$ correspond to a point $\xi\in\bar\fC$ which is
contained in two irreducible components $\bar\fC,\bar\fD$ of $\bar\fX$.
Let $y\in X^\an$ be the point reducing to the generic point of $\bar\fD$
and let $e$ be the edge in $\Sigma(X,V(\fX))$ connecting $x$ and $y$, so
$e$ is a geodesic segment representing $v$.  If $e^\circ$ is the interior
of $e$ then $A = \tau_{V(\fX)}\inv(e^\circ) = \red\inv(\xi)$ in an open annulus;
we let $r$ be the modulus of $A$.
By~\cite[Proposition~3.2]{bosch_lutk:uniformization} we have 
$F(x) - F(y) = -r\cdot\ord_\xi(\td f_x)$.  Since $F$ is affine on $e$ we
also have $F(x) - F(y) = -r \cdot d_v F(x)$, whence the desired
equality.\qed 

\begin{rem}
  Theorem~\ref{thm:PL} 
  is also proved in~\cite[Proposition~3.3.15]{thuillier:thesis},
  in the following form: if $f$ is a nonzero meromorphic
  function on $\hat X^\an$, then the extended real-valued function
  $\log|f|$ on $X$ satisfies the differential equation 
  \begin{equation}
    \label{eq:PL1}
    dd^c \log|f| = \delta_{\Div(f)}
  \end{equation}
  where $dd^c$ is a distribution-valued operator which serves as a
  nonarchimedean analogue of the classical $dd^c$-operator on a Riemann
  surface.  One can regard~\eqref{eq:PL1} as a nonarchimedean analogue of
  the classical `Poincar{\'e}-Lelong formula' for Riemann surfaces.  Since
  it would lead us too far astray to recall the general definition of
  Thuillier's $dd^c$-operator on an analytic curve, we simply call
  Theorem~\ref{thm:PL} the Slope Formula.
\end{rem}

\begin{rem}
  \begin{enumerate}
  \item See~\cite[Example~5.20]{baker_rumely:book} for a version
    of Theorem~\ref{thm:PL} for $X = \bP^1$.
  \item It is an elementary exercise that conditions (4) and (5)
    of Theorem~\ref{thm:PL} uniquely determine  
    the function $F : \Sigma \to \R$ up to addition by a constant; see the proof of
    \cite[Proposition~3.2(A)]{baker_rumely:book}. 
  \end{enumerate}
\end{rem}

\bibliographystyle{thesis}
\bibliography{BPR}
\bigskip~\bigskip

\end{document}